\documentclass[12pt]{amsart}
\usepackage{amsmath,amsthm,amssymb,enumerate,mathscinet,mathtools}
\usepackage{stackengine}
\usepackage{fullpage}
\usepackage{relsize,exscale}
\usepackage[svgnames,dvipsnames]{xcolor}
\usepackage{graphicx,xcolor,pgfplots}
\usepackage{verbatim}
\usepackage{mathrsfs}
\usepackage{hyperref}
\usepackage{enumitem}
\usepackage{blindtext}
\usepackage{multicol}
\usetikzlibrary{fit,calc,positioning,decorations.pathreplacing,matrix}
\newcommand*{\Scale}[2][4]{\scalebox{#1}{$#2$}}
\newcommand{\aaa}{\Scale[1.1]{\mathbf{a}}}
\newcommand{\bbb}{\Scale[1.1]{\mathbf{b}}}
\newcommand{\ccc}{\Scale[1.1]{\mathbf{c}}}
\newcommand{\ddd}{\Scale[1.1]{\mathbf{d}}}
\newcommand{\eee}{\Scale[1.1]{\mathbf{e}}}
\newcommand{\fff}{\Scale[1.1]{\mathbf{f}}}
\newcommand{\gggg}{\Scale[1.1]{\mathbf{g}}}

\newtheorem{theorem}{Theorem}[section]
\newtheorem{fact}[theorem]{Fact}
\newtheorem{lemma}[theorem]{Lemma}

\newtheorem{proposition}[theorem]{Proposition}
\newtheorem{corollary}[theorem]{Corollary}

\newtheorem{thm}[theorem]{Theorem}
\newtheorem{prop}[theorem]{Proposition}

\theoremstyle{definition}
\newtheorem{example}[theorem]{Example}
\newtheorem{remark}[theorem]{Remark}

\numberwithin{equation}{section}

\def\N{\mathbb{N}}

\def\Q{\mathbb{Q}}
\def\W{\mathbb{W}}
\def\R{\mathbb{R}}

\makeatletter
\newcommand*{\rom}[1]{\expandafter\@slowromancap\romannumeral #1@}

\makeatother

\newcommand\scalemath[2]{\scalebox{#1}{\mbox{\ensuremath{\displaystyle #2}}}}

\def\COMMENT#1{}
\let\COMMENT=\footnote% COMMENT OUT for clean output
\allowdisplaybreaks

\title{Enumeration in the lattice of $q$-decreasing words}

\date{\today}

\author{}

\thanks{}
\date{\today}

\begin{document}
\author[J.-L. Baril, N. Hassler, S. Kirgizov]{Jean-Luc Baril, Nathana\"el Hassler, Sergey Kirgizov}
\address{LIB, Universit\'e Bourgogne Europe,   B.P. 47 870, 21078, Dijon Cedex, France}
\email{\{barjl,sergey.kirgizov\}@u-bourgogne.fr, nathanael.hassler@ens-rennes.fr}\theoremstyle{definition}

\begin{abstract}
  We prove that the poset of $q$-decreasing words equipped with the componentwise order forms a lattice. We enumerate the join-irreducible elements for arbitrary $q>0$, and for any positive rational number $q$, we determine the number of coverings, intervals and meet-irreducible elements. The latter present the same structure as words over an alphabet  of $2\lceil q\rceil+1$ letters avoiding $\lceil q\rceil^2+2\lceil q\rceil-1$ consecutive patterns of length 2. Furthermore, we analyze the asymptotic behavior of several of these quantities.
\end{abstract}

\maketitle

\section{Introduction}

Let $q$ be a non-negative real number. A $q$-{\it decreasing word} of length $n\geq 0$ is a binary word of length $n$ satisfying the following constraint: 

\centerline{\it every maximal factor of the form $0^a1^b$ satisfies either $a=0$ or $q\cdot a>b$.} 
More precisely, whenever a block of zeros is followed by a block of ones, the length of the ones block must be strictly less than $q$ times the length of the preceding zeros block (see \cite{Barc,arXivBKV,bkv gray code,egecioglu irsic}). Therefore, such words can start with arbitrarily many 1's, and end with arbitrarily many 0's. Let $\mathcal{W}_n^q$ be the set of $q$-decreasing words of length $n$. For instance, we have $\mathcal{W}_4^1=\{0000,0001,0010,1000,1001,1100,1110,1111\}$, $\mathcal{W}_4^{\frac{1}{2}}=\{0000,0001,1000,1100,1110,1111\}$ and $\mathcal{W}_3^{\frac{\pi}{2}}=\{000,001,010,100,101,110,111\}$. It directly follows  that if $q<r$, then $\mathcal{W}_n^q\subseteq\mathcal{W}_n^r$ for any $n\geq 1$. Notice that if we define $\mathcal{W}_n^{q^+}=\bigcap_{r>q}\mathcal{W}_n^r$, or equivalently, $\mathcal{W}_n^{q^+}$ is the set of binary words of length $n$ such that \\
\centerline{\it every maximal factor of the form $0^a1^b$ satisfies either $a=0$ or $q\cdot a\geq b$,} 
then we have  $\mathcal{W}_n^{q}=\mathcal{W}_n^{q^+}$ for any $n$ when $q$ is irrational, and $\mathcal{W}_n^{q}\varsubsetneq\mathcal{W}_n^{q^+}$ when $q=c/d$, with $c/d$ an irreducible fraction and $n$ sufficiently large ($n\geq c+d$). For instance, we have $$\mathcal{W}_4^{1^+}=\{0000,0001,0010,0100,1000,1001,1010,1100,0101,0011,1110,1101,1111\},$$ which strictly contains $\mathcal{W}_4^1$.

Any word $w\in\mathcal{W}_n^q$ can be written 
        $$w=1^m0^{a_1}1^{b_1}\ldots0^{a_k}1^{b_k}0^\ell,$$
    with $m,\ell\geq 0$, and $q\cdot a_i>b_i\geq1$ for $1\leq i\leq k$ (with 
$k$ possibly equal to zero). Maximal factors  of the form $0^a1^b$ with $q\cdot a>b\geq 1$ will be called {\it prime factors}.
    
The notion of $q$-decreasing words has recently attracted significant attention in the literature. This family of words exhibits a striking combinatorial property whenever $q$ is a positive integer. Indeed, they are in one-to-one correspondence with binary strings that avoid the pattern $1^{q+1}$, i.e. binary strings without $q+1$ consecutive 1 (see \cite{bkv gray code}). So, this implies that $q$-decreasing words of length $n$ (when $q$ is a positive integer) are enumerated by the $(q+1)$-generalized Fibonacci numbers $F_{n+1}^{q+1}$ where $F_n^q$ is defined by 

$$F_n^{q}=F_{n-1}^{q}+F_{n-1}^{q}+\ldots +F_{n-q}^{q},$$ with initial conditions $F_n^q=0$ for $n<0$ and $F_0^q=1$ (see \cite{Fein,knuth,Mile}). It is well known that the generating function of these numbers is $$F_q(x)=\sum_{n\geq 0}F_n^qx^n=\frac{1}{1-x-x^2-\cdots -x^q}.$$

Recently, Barcucci, Bernini, Bilotta and Pinzani \cite{Barc} extended this bijection to $q$-decreasing words for any positive rational number $q$, showing that $\mathcal{W}^q$ is in one-to-one correspondence with binary words avoiding some patterns. 

These words have also been studied from a generative prospective. Baril et al. \cite{bkv gray code} provide efficient algorithms for the generation of all $q$-decreasing words whenever $q$ is a positive integer. In particular, they construct a $3$‑Gray code for general $q$, and notably a $1$‑Gray code for the case $q=1$, thus resolving a conjecture posed in the context of interconnection networks by \cite{egecioglu irsic}. More recently, Wong et al. \cite{Wong} present a two-stage algorithm for generating cyclic $2$-Gray codes for $q$-decreasing words.

More generally, for any $q>0$, the generating function $W_q(x)$ for the number of $q$-decreasing words with respect to the length $n$ is given by 
\begin{equation}\label{eq:w_q general}
W_q(x)=\frac{1}{(1-x)\left(1-\sum_{i=0}^{+\infty}x^{1+i+\left\lfloor\frac{i}{q}\right\rfloor}\right)},    
\end{equation}
see \cite{sergey2,Serg}.
This expression can be simplified as follows when $q$ is rational, i.e. $q=c/d$ where $c$ and $d$ are positive integers:
\begin{equation}\label{eq:w_q rational}
W_q(x)=\frac{1-x^{c+d}}{(1-x)\left(1-x^{c+d}-\sum_{i=0}^{c-1}x^{1+i+\left\lfloor\frac{i}{q}\right\rfloor}\right)}.    
\end{equation}
Note that when $q$ is an integer, i.e. when we fix $d=1$ and $c=q$ in the previous formula,  we obtain $W_q(x)=\frac{F_{q+1}(x)-1}{x}$.  Using (\ref{eq:w_q general}) and (\ref{eq:w_q rational}), Dovgal and Kirgizov \cite{sergey2} proved that for all real $q>0$, $[x^n]W_q(x)\underset{n\to\infty}{\sim}C_q\cdot\Phi(q)^n$, for a positive constant $C_q$, and a function $\Phi(q)$ that interpolates the $q$-bonacci numbers. In particular, when $q=c/d$ is a rational number, $\Phi(q)^{-1}$ is the smallest root in modulus of the polynomial $x^{c+d}+\sum_{i=0}^{c-1}x^{1+i+\left\lfloor\frac{i}{q}\right\rfloor}-1$. See \cite{sergey2} for additional properties of $\Phi(q)$.

To conclude this set of definitions, we introduce the main order-theoretic concepts used throughout this paper. These notions are standard and can be found, for instance, in \cite{Grat,stanley}. A {\it poset} $\mathcal{L}$ is a set endowed with a partial order
relation $\leq$. Given two elements $P,Q\in \mathcal{L}$, a {\it meet} (or
{\it greatest lower bound}) of $P$ and $Q$, denoted $P\wedge Q$, is an
element $R$ such that $R\leq P$, $R\leq Q$, and for any $S$ such that
$S\leq P$ and $S\leq Q$,  then we have $S\leq R$. Dually, a {\it join} (or
{\it least upper bound}) of $P$ and $Q$, denoted $P\vee Q$, is an
element $R$ such that $P\leq R$, $Q\leq R$, and for any $S$ such that
$P\leq S$ and $Q\leq S$, then we have $R\leq S$. Notice that join and meet
elements do not necessarily exist in a poset. A {\it lattice} is a
poset where any pair of elements admits a meet and a join.
An element $P\in\mathcal{L}$ is {\it join-irreducible} (resp. {\it
  meet-irreducible}) if $P=R\vee S$ (resp. $P=R\wedge S$) implies
$P=R$ or $P=S$. 
  An {\it interval} $I$ in a poset $\mathcal{L}$ is a subset of
  $\mathcal{L}$ such that there exist $P,Q\in I$, $P\leq Q$, such that $I=\{R\in\mathcal{L}\ | \ P\leq R \mbox{ and } R\leq Q\}$. An  element $w$ is said to {\it cover} an element $v$  if $v<w$ and there is no $u\in\mathcal{L}$ such that $v<u<w$. In this case, we write  $v\lessdot w$ and the relation between $v$ and $w$ is called a {\it covering}. % A {\it linear interval} $L$ is an interval that is isomorphic to a chain (i.e., if  $x,y\in L$, then $x\leq y$ or $y\leq x$).

In this paper, we endow $\mathcal{W}_n^q$ with the componentwise order. If $v=v_1v_2\ldots v_n$ and  $w=w_1w_2\ldots w_n$ are two $q$-decreasing words in $\mathcal{W}_n^q$ then 
$$v\leq w\Longleftrightarrow v_i\leq w_i \mbox{ for all } 1\leq i\leq n.$$ 

% More formally, if $v=v_1\ldots v_n$ and $w=w_1\ldots w_n$, $v\leq w$ if 
%$$\{i\in [1,n] \ | \ v_i=1\}\subseteq\{i\in [1,n] \ | \ w_i=1\}.$$

\begin{figure}[h]
    \centering
%     \begin{tikzpicture}[scale=1]
%     \node (7) at (0,6) {$1111$};
%     \node (6) at (0,4.5) {$1110$};
%     \node (5) at (2,3) {$1001$};
%     \node (4) at (0,3) {$1100$};
%     \node (3) at (2,1.5) {$0001$};
%     \node (2) at (0,1.5) {$1000$};
%     \node (1) at (-2,1.5) {$0010$};
%     \node (0) at (0,0) {$0000$};
%     \draw (0) -- (1);\draw (0) -- (2);\draw (0) -- (3);\draw (1) -- (6);\draw (2) -- (4);\draw (2) -- (5);\draw (3) -- (5);\draw (4) -- (6);\draw (5) -- (7);\draw (6) -- (7);
% \end{tikzpicture}
\begin{tikzpicture}[scale=1.3]
    \node (11111) at (0,7.5) {$11111$};
    \node (11110) at (-1.5,6) {$11110$};
    \node (11001) at (0,4.5) {$11001$};
    \node (11100) at (-3,4.5) {$11100$};
    \node (00011) at (3,3) {$00011$};
    \node (10001) at (1,3) {$10001$};
    \node (10010) at (-1,3) {$10010$};
    \node (11000) at (-3,3) {$11000$};
    \node (00001) at (3,1.5) {$00001$};
    \node (00010) at (1,1.5) {$00010$};
    \node (10000) at (-3,1.5) {$10000$};
    \node (00100) at (-1,1.5) {$00100$};
    \node (00000) at (0,0) {$00000$};
    \draw (00000) -- (00100);\draw (00000) -- (10000);\draw (00000) -- (00001);\draw (00000) -- (00010);\draw (10000) -- (11000);\draw (10000) -- (10010);\draw (10000) -- (10001);\draw (00010) -- (10010);\draw (00010) -- (00011);\draw (00001) -- (10001);\draw (00001) -- (00011);\draw (00100) -- (11100);\draw (11000) -- (11001);\draw (11000) -- (11100);\draw (10010) -- (11110);\draw (10001) -- (11001);\draw (00011) -- (11111);\draw (11100) -- (11110);\draw (11001) -- (11111);\draw (11110) -- (11111);
\end{tikzpicture}
    \caption{The lattice $\W_5^1$. It contains 20 coverings (edges), 7 meet-irreducible elements, 5 join-irreducible elements and 56 intervals. }
    \label{fig:W_4^2}
\end{figure}
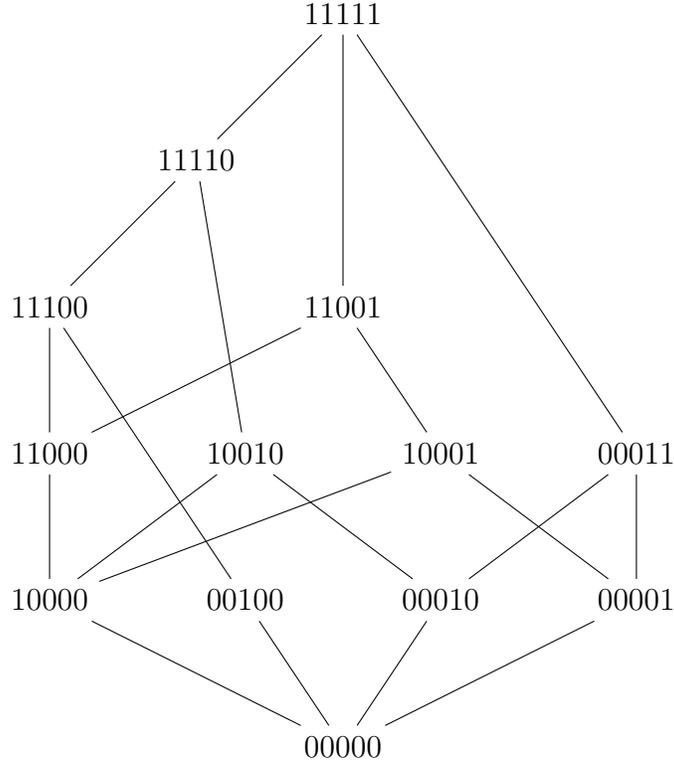

Let  $\W_n^q:=(\mathcal{W}_n^q,\leq)$ be the poset defined by this order relation.  See Figure~\ref{fig:W_4^2} for an illustration of the poset $\W_5^1$.

\begin{thm} For $q\geq 0$, the poset     $\W_n^q$ is a lattice for any $n\geq 1$.
\end{thm}
\begin{proof} For any words $v$ and $w$ in $\mathcal{W}_n^q$, we consider the binary word $a=a_1\ldots a_n$  where $a_i=1$ if and only if $v_i=w_i=1$. It is straightforward to see that $a\in \mathcal{W}_n^q$, and then $a$ is the greatest lower bound of $v$ and $w$,  which implies that $\W_n^q$ is a meet-semilattice. Since $1^n$ is the maximum element of $\W_n^q$, Proposition 3.3.1 in \cite{stanley} implies that $\W_n^q$ is a lattice.   
\end{proof}

\noindent {\bf Outline of the paper.} In Section~\ref{sec:useful results}, we collect preliminary results that will be used throughout the paper. Many of these results are quite technical, due to the presence of floor and ceiling functions in generating functions related to $q$-decreasing words. Section~\ref{sec:join and coverings} is devoted to enumerative results concerning classical lattice parameters, specifically join-irreducible elements, as well as covering relations. For any rational number \( q > 0 \), we provide a closed-form expression for the generating function counting the number of coverings in \( \W^q_n \). When \( q \) is irrational, we present a formula that enables efficient computation of the initial terms of the series expansion (e.g., using \textsc{Maple}). We also prove that the asymptotic behavior of the number of coverings is connected to the function $\Phi(q)$ defined in the introduction above. In Section~\ref{sec:intervals}, we derive a closed-form expression for the generating function that enumerates the number of intervals in \( \W^q_n \) for any rational number \( q > 0 \). Finally, Section~\ref{sec:meet} presents the structure of meet-irreducible elements in $\W_n^q$ for any positive rational number $q$. This structure is the same as the one of words over an alphabet  of $2\lceil q\rceil+1$ letters avoiding $\lceil q\rceil^2+2\lceil q\rceil-1$ consecutive patterns of length 2. Taking advantage of this classical structure, we present a method to obtain the closed form of the generating function enumerating the number of meet-irreducible elements in $\W_n^q$ for any rational number $q>0$.

\section{Useful results}\label{sec:useful results}

In this section, we collect various results that will help us in the study of the generating functions we obtain in Section \ref{sec:join and coverings}, \ref{sec:intervals} and \ref{sec:meet}. We begin with two simple facts about floor and ceiling functions.

\begin{fact}\label{fact partie entière}
Let $a,b$ be integers with $a\geq 0$ and $b>0$, then
    $$1+\left\lfloor\frac{a}{b}\right\rfloor=\left\lceil\frac{a+1}{b}\right\rceil.$$
\end{fact}
\begin{proof}
    It follows from the following inequalities $\left\lfloor\frac{a}{b}\right\rfloor\leq\frac{a}{b}<\frac{a+1}{b}\leq \left\lfloor\frac{a}{b}\right\rfloor+1$. 
\end{proof}
% \begin{proof}
%     Let $n=kc+r$ with $k\geq0$ and $0\leq r < c$ be the euclidian division of $n$ by $c$. Then
%     $$\frac{(c+d)n}{c}=(c+d)k+\frac{(c+d)r}{c}, \mbox{ and } \frac{(c+d)n+1}{c}=(c+d)k+\frac{(c+d)r+1}{c}.$$
%     So it suffices to check that 
%     $1+\left\lfloor\frac{(c+d)r}{c}\right\rfloor=\left\lceil\frac{(c+d)r+1}{c}\right\rceil$, which is equivalent to $1+\left\lfloor\frac{dr}{c}\right\rfloor=\left\lceil\frac{dr+1}{c}\right\rceil$. This is true if $r=0$, and otherwise, since $c$ and $d$ are relatively prime and $0<r<c$, $c$ does not divide $dr$. This means that $\frac{dr}{c}$ is not an integer, so there exists an integer $m$ such that $m<\frac{dr}{c}<m+1$, thus $m<\frac{dr+1}{c}\leq m+1$. Then, $m=\left\lfloor\frac{dr}{c}\right\rfloor$ and $m+1=\left\lceil\frac{dr+1}{c}\right\rceil$.
% \end{proof}
The following fact characterizes the integers $b$ such that $\left\lfloor\frac{b}{q}\right\rfloor=\left\lfloor\frac{b+1}{q}\right\rfloor$.
\begin{fact}\label{carac b/q=(b+1)/q}
    Let $q=c/d$ and $b\geq1$. Then
    $$\left\lfloor\frac{b}{q}\right\rfloor=\left\lfloor\frac{b+1}{q}\right\rfloor\Leftrightarrow db\bmod{c}\in\{0,\ldots,c-d-1\}.$$
\end{fact}
\begin{proof}
Let $db=sc+r$, with $s\geq0$ and $0\leq r<c$ be the euclidean division of $db$ by $c$. Then $\frac{b}{q}=s+\frac{r}{c}$, and $\frac{b+1}{q}=s+\frac{r+d}{c}$. So, we have
$$s=\left\lfloor\frac{b}{q}\right\rfloor=\left\lfloor\frac{b+1}{q}\right\rfloor\Leftrightarrow \frac{r+d}{c}<1\Leftrightarrow 0\leq r<c-d.$$
\end{proof}
The following technical lemma will help us to characterize the prime factors of meet-irreducible elements in Theorem \ref{mi 0<q<=1}.
\begin{lemma}\label{factor mi each length}
        Let $q=c/d>0$ be a rational number. For every $n\geq2+\left\lfloor\frac{1}{q}\right\rfloor$, there exists a unique pair $(a,b)$ of positive integers  such that $a+b=n$ and $1+\left\lfloor\frac{b}{q}\right\rfloor\leq a\leq 1+\left\lfloor\frac{b+1}{q}\right\rfloor$.
    \end{lemma}
\begin{proof}
Let $b=\left\lceil\frac{qn}{q+1}\right\rceil-1$, and $a=n-b$. Since $n\geq2+\left\lfloor\frac{1}{q}\right\rfloor$, we have $1\leq a,b\leq n-1$. By Fact \ref{fact partie entière}, $b=\left\lceil\frac{cn}{c+d}\right\rceil-1=\left\lfloor\frac{cn-1}{c+d}\right\rfloor$. Then we have
\begin{align*}
    1+\left\lfloor\frac{b}{q}\right\rfloor
    \leq 1+\frac{db}{c}=1+\frac{d}{c}\left\lfloor\frac{cn-1}{c+d}\right\rfloor
    \leq 1+\frac{d}{c}\cdot\frac{cn-1}{c+d}
    =\frac{dn}{c+d}+1-\frac{d}{c(c+d)}.
\end{align*}
Since $1+\left\lfloor\frac{b}{q}\right\rfloor$ is an integer, we actually have 
$$1+\left\lfloor\frac{b}{q}\right\rfloor\leq\left\lfloor\frac{dn}{c+d}+1-\frac{d}{c(c+d)}\right\rfloor\leq 1+\left\lfloor\frac{dn}{c+d}\right\rfloor=\left\lceil\frac{dn+1}{c+d}\right\rceil,$$
where we used Fact \ref{fact partie entière} for the last equality. On the other hand, 
\begin{align*}
    a=n-b=n-\left\lfloor\frac{cn-1}{c+d}\right\rfloor=\left\lceil n-\frac{cn-1}{c+d}\right\rceil=\left\lceil\frac{dn+1}{c+d}\right\rceil\geq1+\left\lfloor\frac{b}{q}\right\rfloor.
\end{align*}
For the second inequality, 
\begin{align*}
    1+\left\lfloor\frac{b+1}{q}\right\rfloor&=1+\left\lfloor\frac{db+d}{c}\right\rfloor=\left\lceil\frac{db+d+1}{c}\right\rceil\geq\frac{db+d+1}{c}=\frac{d+1}{c}+\frac{d}{c}\left(\left\lceil\frac{cn}{c+d}\right\rceil-1\right)\\
    &\geq\frac{d+1}{c}+\frac{d}{c}\cdot\frac{cn}{c+d}-\frac{d}{c}=\frac{dn}{c+d}+\frac{1}{c}=\frac{cdn+c+d}{c(c+d)}.
\end{align*}
Since $1+\left\lfloor\frac{b+1}{q}\right\rfloor$ is an integer and by Fact \ref{fact partie entière}, we actually have 
$$1+\left\lfloor\frac{b+1}{q}\right\rfloor\geq \left\lceil\frac{cdn+c+d}{c(c+d)}\right\rceil=1+\left\lfloor\frac{cdn+c+d-1}{c(c+d)}\right\rfloor=1+\left\lfloor\frac{dn}{c+d}+\frac{1}{c}-\frac{1}{c(c+d)}\right\rfloor\geq 1+\left\lfloor\frac{dn}{c+d}\right\rfloor.$$
On the other hand, 
$$a=n-b=n-\left\lceil\frac{cn}{c+d}\right\rceil+1=\left\lfloor n-\frac{cn}{c+d}\right\rfloor+1=\left\lfloor\frac{dn}{c+d}\right\rfloor+1\leq 1+\left\lfloor\frac{b+1}{q}\right\rfloor.$$
% For the uniqueness, we will see in Theorem \ref{mi 0<q<=1} that 
% $$1+\left\lfloor\frac{b}{q}\right\rfloor\leq a\leq 1+\left\lfloor\frac{b+1}{q}\right\rfloor\Leftrightarrow 0^a1^b \mbox{ is a factor of a meet-irreducible element in } \W_{a+b}^q.$$
% Now suppose that $w=1^m0^{a_1}1^{b_1}\ldots0^a1^b\ldots0^{a_k}1^{b_k}0^\ell$ and $w'=1^m0^{a_1}1^{b_1}\ldots0^{a'}1^{b'}\ldots0^{a_k}1^{b_k}0^\ell$ are both meet-irreducible, with $a+b=a'+b'=n$, and $a<a'$. Then we would have $w'\leq w$ and $w'\leq 1^{m+a_1+b_1}\ldots0^{a'}1^{b'}\ldots0^{a_k}1^{b_k}0^\ell$, contradicting the meet-irreducibility of $w'$, since $w$ and $1^{m+a_1+b_1}\ldots0^{a'}1^{b'}\ldots0^{a_k}1^{b_k}0^\ell$ are not comparable.

For uniqueness, suppose that $(a,b)$, $(a',b')$ are two such pairs with $a<a'$, and thus $b>b'$. Then we would have
$$1+\left\lfloor\frac{b'+1}{q}\right\rfloor\leq1+\left\lfloor\frac{b}{q}\right\rfloor\leq a\leq a'-1\leq \left\lfloor\frac{b'+1}{q}\right\rfloor,$$
which is a contradiction.
\end{proof}

In the following lemma, we determine the generating function of the sequence $(\lfloor \frac{an+c}{b}\rfloor)_{n\in\N}$ where $a,b$ and $c$ are  some integers satisfying $a,b>0$, and we show that it can be expressed as a rational function.

\begin{lemma}\label{gf (an/b)}
    Let $a,b$ be two positive integers, relatively prime, and let $c$ be an integer. Then we have 
    $$\sum_{n=0}^{+\infty}\left\lfloor\frac{an+c}{b}\right\rfloor x^n=\frac{ax^b+(1-x)Q_{a,b,c}(x)}{(1-x)(1-x^b)},$$
    with $Q_{a,b,c}(x)=\sum_{r=0}^{b-1}\left\lfloor\frac{ar+c}{b}\right\rfloor x^r$.
\end{lemma}
\begin{proof}
    Let $n\in\N$, and let $n=bk+r$ be its euclidean division by $b$. Thus we have $0\leq r<b$. A simple calculation provides 
    $$\left\lfloor\frac{an+c}{b}\right\rfloor=ak+\left\lfloor\frac{ar+c}{b}\right\rfloor,$$ and then
    \begin{align*}
    \sum_{n=0}^{+\infty}\left\lfloor\frac{an+c}{b}\right\rfloor x^n&=\sum_{k=0}^{+\infty}\sum_{r=0}^{b-1}\left(ak+\left\lfloor\frac{ar+c}{b}\right\rfloor\right)x^{bk+r}\\
    &=\sum_{k=0}^{+\infty}x^{bk}\left(ak\sum_{r=0}^{b-1}x^r+Q_{a,b,c}(x)\right)\\
    &=a\frac{1-x^b}{1-x}\sum_{k=0}^{+\infty}k(x^{b})^k+Q_{a,b,c}(x)\sum_{k=0}^{+\infty}(x^{b})^k\\
    &=\frac{ax^b}{(1-x)(1-x^b)}+\frac{Q_{a,b,c}(x)}{1-x^b}.
    \end{align*}
\end{proof}
\begin{example}
For instance, if $a=3$, $b=5$ and $c=2$, then we have $Q_{3,5,2}(x)=x+x^2+2x^3+2x^4$, and 
$$\sum_{n=0}^{+\infty}\left\lfloor\frac{3n+2}{5}\right\rfloor x^n=\frac{3x^5+(1-x)(x+x^2+2x^3+2x^4)}{(1-x)(1-x^5)}={\frac {x \left(x^4+ {x}^{2}+1 \right) }{ \left( 1-x \right)  \left( 1-{x}^
{5} \right) }}
.$$
\end{example}

The following lemma plays a key role in the computation of generating functions involving a ceiling function in the exponent of $x$. 

% \begin{remark}
%     $n=b+\left\lceil\frac{b+1}{q}\right\rceil\Longrightarrow b=\left\lfloor\frac{q(n+1)}{q+1}\right\rfloor-1$.
% \end{remark}

\begin{lemma}\label{lemma periodicity}
Let $g(x)=\sum_{n=0}^{+\infty} a_n x^n$, and let $q=\frac{c}{d}$ be a positive rational number. We define the series $g_{c/d}$ as $$g_{c/d}(x)=\sum_{n=0}^{+\infty}a_nx^{1+\left\lfloor\frac{(c+d)n}{c}\right\rfloor}.$$ Then we have 
$$g_{c/d}(x)=\frac{1}{c}\sum_{j=0}^{c-1}\sum_{k=0}^{c-1}\rho^{-kj}g(\rho^kx^{\frac{c+d}{c}})x^{1+\left\lfloor\frac{jd}{c}\right\rfloor-\frac{jd}{c}},$$
where $\rho=e^{2i\pi/c}$ is a primitive root of unity.
\end{lemma}
\begin{proof}
We start with the following two simple facts. For $n\geq 0$, we have 
$$\frac{1}{c}\sum_{k=0}^{c-1} \rho^{kn}=\left\{
\begin{array}{cl}
    1 & \mbox{ if } n\equiv 0 \ [c]\\
    0 & \mbox{ otherwise.}
\end{array}\right.$$
Furthermore, for $j\in [0,c-1]$, we have
\begin{equation}\label{eq1}
    \frac{1}{c}\sum_{k=0}^{c-1}\rho^{-kj}g(\rho^kx)=\sum_{n=0}^{+\infty}\left(\frac{1}{c}\sum_{k=0}^{c-1}\rho^{k(n-j)}\right)a_nx^n=\sum_{n=0}^{+\infty}a_{cn+j}x^{cn+j}.
\end{equation}
Let $n\geq 0$, and let $n=cm+j$ be its euclidean division by $c$. Then, $$1+\left\lfloor\frac{(c+d)n}{c}\right\rfloor=1+(c+d)m+\left\lfloor\frac{(c+d)j}{c}\right\rfloor.$$
Thus, the series $\sum_{n=0}^{+\infty}a_nx^{1+\left\lfloor\frac{(c+d)n}{c}\right\rfloor}$ can be rewritten 
$$\sum_{n=0}^{+\infty}a_nx^{1+\left\lfloor\frac{(c+d)n}{c}\right\rfloor}=\sum_{m=0}^{+\infty}\sum_{j=0}^{c-1}a_{cm+j}x^{1+(c+d)m+\left\lfloor\frac{(c+d)j}{c}\right\rfloor}.$$
It then follows from (\ref{eq1}) that
$$
\sum_{n=0}^{+\infty}a_nx^{1+\left\lfloor\frac{(c+d)n}{c}\right\rfloor}=\sum_{j=0}^{c-1}\sum_{m=0}^{+\infty}a_{cm+j}x^{1+(c+d)m+\left\lfloor\frac{(c+d)j}{c}\right\rfloor}=\sum_{j=0}^{c-1}\frac{x^{1+\left\lfloor\frac{jd}{c}\right\rfloor-\frac{jd}{c}}}{c}\sum_{k=0}^{c-1}\rho^{-kj}g(\rho^kx^{\frac{c+d}{c}}).$$
\end{proof}

In the next proposition, we examine the special case where $g$ is a rational function. In particular, we provide a more efficient method for computing  $g_{c/d}$ than the double sum given in  Lemma \ref{lemma periodicity}. Moreover we derive the  asymptotic behavior of the coefficients of $g_{c/d}$ under the assumption that $g$ has a unique singularity of minimal modulus.

\begin{proposition}\label{closed form transfo}
    Let $g(x)=\frac{P(x)}{Q(x)}$ be a rational power series, and let $c/d$ be a positive rational number, with $c$ and $d$ relatively prime. Then the following two statements hold.
    \begin{enumerate}
        \item $g_{c/d}$ is rational, and $g_{c/d}(x)=N(x)/D(x)$, where
        $$\left\{
        \begin{array}{ll}
        D(x)=\prod_{k=0}^{c-1}Q(\rho^kx^{1+d/c}),&\quad \mbox{ with } \rho=e^{\frac{2i\pi}{c}}, \\
        N(x)=\sum_{k=0}^m a_kx^{1+\left\lfloor\frac{(c+d)k}{c}\right\rfloor},&\quad  \mbox{ where } P(x)\prod_{k=1}^{c-1}Q(\rho^kx)=\sum_{k=0}^ma_kx^k.
        \end{array}
        \right.$$
        \item Suppose further that $[x^n]g(x)\underset{n\to\infty}{\sim}C\cdot\alpha^{-n}$ for some constant $C$ and $\alpha>0$. Then there exists a nonzero polynomial $T\in\mathbb{C}[X]$ with degree at most $c+d-1$ such that
    $$[x^n]g_{c/d}(x)\underset{n\to\infty}{\sim}T(r^n)\cdot\alpha^{-\frac{cn}{c+d}},$$
    with $r=e^{\frac{2i\pi}{c+d}}.$
    \end{enumerate}
\end{proposition}
\begin{proof}
    (1) Let $g(x)=\frac{P(x)}{Q(x)}$ with $P$ and $Q$ relatively prime. Let $\rho=e^{\frac{2i\pi}{c}}$, and for $k=~0,\ldots,c-~1$, let $R_k(x)=\sum_{j=0}^{c-1}\rho^{-kj}x^{1+\left\lfloor\frac{jd}{c}\right\rfloor-\frac{jd}{c}}$ (this is a polynomial in $x^\frac{1}{c}$). By Lemma \ref{lemma periodicity}, 
    \begin{align*}
        g_{c/d}(x)=\frac{1}{c}\sum_{i=0}^{c-1}R_k(x)\frac{P(\rho^kx^{1+d/c})}{Q(\rho^kx^{1+d/c})}=\frac{N(x)}{D(x)},
    \end{align*}
    where $D(x)=\prod_{k=0}^{c-1}Q(\rho^kx^{1+d/c})$, and 
    $$N(x)=\frac{1}{c}\sum_{k=0}^{c-1}R_k(x)P(\rho^kx^{1+d/c})\prod_{j\ne k}Q(\rho^jx^{1+d/c}).$$
    We begin by analyzing $D(x)$. To this end, suppose first that $Q(x)=x-\alpha$. Then we have
    \begin{align*}
        D(x)&=(x^{1+d/c}-\alpha)(\rho x^{1+d/c}-\alpha)\ldots(\rho^{c-1} x^{1+d/c}-\alpha)\\
        &=\left(\prod_{k=0}^{c-1}\rho^k\right)x^{c+d}-\alpha\left(\sum_{k=0}^{c-1}\prod_{\ell\ne k}\rho^\ell\right)x^{\frac{(c-1)(d+c)}{c}}
        +\ldots\\
        &+(-\alpha^j)\left(\sum_{0\leq k_1<\ldots <k_j\leq c-1}\prod_{\ell\not\in\{k_1,\ldots,k_j\}}\rho^\ell\right)x^{\frac{(c-j)(d+c)}{c}}+\ldots+(-\alpha)^c.
    \end{align*}
    It follows from the Vieta's formulas for the polynomial $x^c-1$ that the above expression simplifies to $D(x)=(-1)^{c+1}(x^{c+d}-\alpha^c)$. Thus, for a general polynomial $Q(x)=\lambda\prod_{k=1}^s(x-~\alpha_k)^{m_k}$, we have
    \begin{equation}\label{eq:denominateur equivalent}
        D(x)=\lambda^c\prod_{k=1}^s(-1)^{(c+1)m_k}(x^{c+d}-\alpha_k^c)^{m_k}.
    \end{equation}
    In particular, $D(x)$ is a polynomial. Now let us focus on $N(x)$. Let $S(x)=P(x)\prod_{k=1}^{c-1}Q(\rho^kx)$. Then we have 
    $$N(x)=\frac{1}{c}\sum_{k=0}^{c-1}R_k(x)S(\rho^kx^{1+d/c}).$$
    Therefore, by Lemma \ref{lemma periodicity}, if $S(x)=\sum_{k=0}^m a_kx^k$, then $N(x)=\sum_{k=0}^m a_kx^{1+\left\lfloor\frac{(c+d)k}{c}\right\rfloor}$, which is also a polynomial.\\
    % Moreover, since $g$ is a power series, $\frac{1}{D(x)}$ is a power series too, because none of the $\alpha_k$'s is zero. Since $g_{c/d}$ is also a power series by definition, we deduce that $N(x)=g_{c/d}(x)D(x)$ is a power series. However this is a polynomial in $x^{\frac{1}{c}}$, so $N(x)$ is actually a polynomial. This proves that $g_{c/d}$ is rational.\\
    (2) Now suppose that $[x^n]g(x)\underset{n\to\infty}{\sim}C\cdot\alpha^{-n}$ for some constants $C$ and $\alpha>0$. Then $\alpha$ is the unique root of $Q(x)$ of minimal modulus. By (\ref{eq:denominateur equivalent}), the smallest roots of $D(x)$ are then $\alpha^{\frac{c}{c+d}}, r\alpha^{\frac{c}{c+d}},\ldots,r^{c+d-1}\alpha^{\frac{c}{c+d}}$, with $r=e^{\frac{2i\pi}{c+d}}$. We can check that at least $\alpha^{\frac{c}{c+d}}$ is a singularity of $g_{c/d}$, because
    \begin{align*}
        N(\alpha^{\frac{c}{c+d}})&=\frac{1}{c}R_0(\alpha^{\frac{c}{c+d}})P(\alpha)\prod_{j\ne 0}Q(\rho^j\alpha)\ne 0.
    \end{align*}
    Indeed, for $j\ne0$, we have $Q(\rho^j\alpha)\ne0$ since $\alpha$ is the unique root of $Q$ with modulus $|\alpha|$; moreover, $P(\alpha)\ne0$ because $P$ and $Q$ are relatively prime; and finally, $R_0(x)=\sum_{j=0}^{c-1}\alpha^{\frac{c}{c+d}\left(1+\left\lfloor\frac{jd}{c}\right\rfloor-\frac{jd}{c}\right)}>0$, since $\alpha>0$. We deduce that there exist some constants $c_0,\ldots,c_{c+d-1}$, with $c_0\ne0$ such that 
    $$[x^n]g_{c/d}(x)\underset{n\to\infty}{\sim} (c_0+c_1r^n+\ldots+c_{c+d-1}r^{(c+d-1)n})\alpha^{-\frac{cn}{c+d}}.$$
\end{proof}

\begin{example}
    Let $g(x)=\frac{1}{1-\alpha x}$. Then, by Lemma \ref{lemma periodicity}, 
    $$g_2(x)=\frac{1}{2}\left(\frac{x+x^{1/2}}{1-\alpha x^{3/2}}+\frac{x-x^{1/2}}{1+\alpha x^{3/2}}\right)=\frac{x(1+\alpha x)}{1-\alpha^2x^3}.$$
    Using routine singularity analysis (see e.g. \cite{flaj}), we deduce that 
    $$[x^n]g_2(x)\underset{n\to\infty}{\sim}(c_0+c_1j^n+c_2j^{2n})\alpha^{2n/3},$$
    where $j=e^{\frac{2i\pi}{3}}$ and $c_k=\frac{j^k\alpha^{-1/3}(1+j^k\alpha^{-1/3})}{3}$ for $k=0,1,2$. Furthermore,
    $$c_0+c_1j^n+c_2j^{2n}=\left\{\begin{array}{cc}
        0 &\mbox{ if } n=0\bmod{3},   \\
        \alpha^{-2/3}  &\mbox{ if } n=1\bmod{3},   \\
        \alpha^{-1/3}  &\mbox{ if } n=2\bmod{3}. 
    \end{array}\right.$$
\end{example}

\section{Join-irreducible elements and coverings}\label{sec:join and coverings}

In this section, we provide enumerative results for the classical parameters of a lattice, namely the join-irreducible elements, as well as covering relations. The enumeration of meet-irreducible elements is a bit more intricate, so we treat it in Section \ref{sec:meet}.
We first give the enumeration of join-irreducible elements for any $q>0$, and then we conclude by giving closed form for the generating functions of the covering for positive rational numbers $q$, and a method for computing arbitrarily many terms of the generating functions for positive irrational numbers $q$.  
\subsection{Join-irreducible elements}
\begin{thm}\label{join irreducible}
    For $q>0$ and $n\geq1$, there are exactly $n$ join-irreducible elements in $\W_n^q$.
\end{thm}
\begin{proof}
    In a finite lattice, an element is join-irreducible if and only if it covers exactly one element. We therefore count the words in $\W_n^q$ that cover exactly one element. Note that a factor $0^a1^b$, with $qa>b\geq 1$ covers only one element if and only if $b=1$ (otherwise it covers $0^a1^{b-1}0$ and $0^{a+1}1^{b-1}$). Now we investigate the elements covered by $1^m$. Given $i\geq0$, there is at most one word covered by $1^m$ with the suffix $01^i$. So the words covered by $1^m$ are exactly

    $$1^{n-1}0,\,1^{n-\alpha_1-1}0^{\alpha_1}1,\ldots,\,1^{n-\alpha_k-k}0^{\alpha_k}1^k,$$
    where $\alpha_i=1+\left\lfloor\frac{i}{q}\right\rfloor$ is the smallest integer such that $q\cdot \alpha_i>i$, and $k=\left\lceil\frac{qm}{q+1}\right\rceil-1$ is the largest integer such that $\alpha_k\leq m-k$. Thus, a factor $1^m$ covers only one element if and only if $\left\lceil\frac{qm}{q+1}\right\rceil=1$, i.e. if and only if $1\leq m\leq 1+\left\lfloor\frac{1}{q}\right\rfloor$. So, taking a word $w\in\W_n^q$, and its decomposition
    $$w=1^m0^{a_1}1^{b_1}\ldots0^{a_k}1^{b_k}0^\ell,$$
    $w$ is join-irreducible if and only if $m=0$, $k=1$, $b_1=1$ and $qa_1>1$, or $1\leq m\leq 1+\left\lfloor\frac{1}{q}\right\rfloor$ and $k=0$. Thus, the generating function for the number of join-irreducible elements in $\W_n^q$ is
    $$\frac{\sum_{qa>1}x^{a+1}}{1-x}+\frac{x+x^2+\ldots+x^{1+\left\lfloor\frac{1}{q}\right\rfloor}}{1-x}=\frac{x^{2+\left\lfloor\frac{1}{q}\right\rfloor}}{(1-x)^2}+\frac{x-x^{2+\left\lfloor\frac{1}{q}\right\rfloor}}{(1-x)^2}=\frac{x}{(1-x)^2}.$$
    This proves that for any $q>0$, there are $n$ join-irreducible elements in $\W_n^q$.
\end{proof}

\subsection{Coverings}

\begin{proposition}\label{structure coverings}
    Let $v,w\in\mathcal{W}_n^q$. Writing  $w=1^m0^{a_1}1^{b_1}\ldots0^{a_k}1^{b_k}0^\ell$   with $m,\ell\geq 0$, and $q\cdot a_i>b_i\geq1$ for $1\leq i\leq k$, we suppose that $v=u_1v_1\ldots v_ku_0,$
    with $|u_1|=m$, $|u_0|=\ell$, and $|v_i|=a_i+b_i$ for $1\leq i\leq k$. 
    Then $w$ covers $v$ ($v\lessdot w$) if and only if $u_0=0^\ell$, and
    one of the following two statements holds:

    \begin{itemize}
        \item[(1)] $u_1\lessdot 1^m \mbox{ and } v_i=0^{a_i}1^{b_i} \mbox{ for } 1\leq i\leq k$,
        \item[(2)] $u_1= 1^m \mbox{ and there exists a unique  } i \mbox{ such that } v_i\lessdot0^{a_i}1^{b_i}, \mbox{ and } v_j=0^{a_j}1^{b_j} \mbox{ for } j\ne i$.
    \end{itemize}
    %$$\left\{\begin{array}{ll}
     %    &\nu_1\lessdot 1^d \mbox{ and } v_i=0^{a_i}1^{b_i} \mbox{ for } 1\leq %i\leq k, \\
     %    \mbox{ or } &  
    %\end{array}\right.$$
\end{proposition}

\begin{proof}
The implication from right to left being easy to check, we focus on the converse. Suppose that $v\lessdot w$, i.e. $w$ covers $v$. In particular, $u_1\leq 1^m$, $u_0\leq 0^\ell$ and $v_i\leq 0^{a_i}1^{b_i}$ for  $1\leq i\leq k$. Then we necessarily have $u_0=0^\ell$. Now assume that among $u_1,v_1,\ldots,v_k$, there are at least two of them (say for instance and without loss of generality $u_1$ and $v_1$) such that $u_1<1^m$ and $v_1< 0^{a_1}1^{b_1}$. Then we would have $v<1^mv_1\ldots v_k 0^\ell<w$, contradicting $v\lessdot w$. If we had chosen $v_1< 0^{a_1}1^{b_1}$  and $v_2< 0^{a_2}1^{b_2}$, then the inequalities $v<1^m0^{a_1}1^{b_1}v_2\ldots v_k 0^\ell<w$ would have also led to a contradiction. So there is exactly one word among $u_1,v_1,\ldots,v_k$ such that $u_1\lessdot 1^m$ or $v_i\lessdot 0^{a_i}1^{b_i}$.
\end{proof}

Let $A_q(x)$ be the generating function for the number of coverings of the form $\nu\lessdot 1^m$. Let $B_q(x)$ denote the generating function for the number of coverings of the form $v\lessdot0^a1^b$ with $q\cdot a>b\geq1$. Let $D_q(x)$ denote the generating function for the number of elements in $\mathcal{W}_n^q$ of the form $0^a1^b$ with $q\cdot a>b\geq1$.

\begin{proposition}\label{nb coverings general}
     Let $q>0$ be a real number, and let $C_q(x)$ denote the generating function for the number of coverings in $\W_n^q$. Then 
    $$C_q(x)=\frac{A_q(x)}{(1-x)(1-D_q(x))}+\frac{B_q(x)}{(1-x)^2(1-D_q(x))^2}.$$
\end{proposition}
\begin{proof}
    Using Proposition \ref{structure coverings},  coverings of the form $(1)$ are in one-to-one correspondence with words $v=u_1v_1v_2\ldots v_k0^\ell$ satisfying $u_1\lessdot 1^m$ and $v_i=0^{a_i}1^{b_i}$ for $1\leq i\leq k$; thus they contribute to  $$\frac{A_q(x)}{(1-x)(1-D_q(x))}.$$ Coverings of the form $(2)$ are in one-to-one correspondence with words $v=1^mv_1v_2\ldots v_k0^\ell$ such that there is a unique $i$ with $v_i\lessdot 0^{a_i}1^{b_i}$ and $v_j=0^{a_j}1^{b_j}$ for $j\neq i$; thus they contribute to $$\frac{B_q(x)}{(1-x)^2(1-D_q(x))^2}.$$ Considering these two cases, the expected result follows.
\end{proof}

The next proposition establishes the explicit forms of $A_q(x)$, $B_q(x)$ and $D_q(x)$ for any $q>0$.

\begin{proposition}\label{nb coverings}
    For any real $q>0$, we have 
    $$A_q(x)=\sum_{n=1}^{+\infty}\left\lceil\frac{qn}{q+1}\right\rceil x^n,\quad B_q(x)=\frac{1}{1-x}\sum_{b=1}^{+\infty}\left\lceil\frac{qb+1}{q+1}\right\rceil x^{1+b+\left\lfloor\frac{b}{q}\right\rfloor},$$
    
    $$\mbox{and}\quad D_q(x)=\frac{1}{1-x}\sum_{b=1}^{+\infty}x^{1+b+\left\lfloor\frac{b}{q}\right\rfloor}.$$
\end{proposition}
\begin{proof} The proof is divided into three parts, each corresponding to $A_q(x)$, $B_q(x)$ and $D_q(x)$ respectively. 

\begin{itemize}
    \item As we saw in the proof of Theorem \ref{join irreducible}, for $n\geq 1$, the words $v\in\mathcal{W}_n^q$ such that $v\lessdot 1^n$ are exactly $$1^{n-1}0,\,1^{n-\alpha_1-1}0^{\alpha_1}1,\ldots,\,1^{n-\alpha_k-k}0^{\alpha_k}1^k,$$
    where $\alpha_i=1+\left\lfloor\frac{i}{q}\right\rfloor$ is the smallest integer such that $q\cdot \alpha_i>i$, and $k=\left\lceil\frac{qn}{q+1}\right\rceil-1$ is the largest integer such that $\alpha_k\leq n-k$. Indeed for each $i\in[0,k]$, there is exactly one word covered by $1^n$ with suffix $01^i$. So there are exactly
    $\left\lceil\frac{qn}{q+1}\right\rceil$
    words covered by $1^n$ in $\W_n^q$, hence the expression of $A_q(x)$.
    
    \item Let $a,b$ be positive integers such that $q\cdot a>b\geq1$. The words $v\in\mathcal{W}_n^q$ such that $v\lessdot 0^a1^b$ are exactly
    $$0^a1^{b-1}0,\,0^a1^{b-\alpha_1-1}0^{\alpha_1}1,\ldots,\,0^a1^{b-\alpha_k-k}0^{\alpha_k}1^k,\,0^{a+1}1^{b-1},$$
    where $\alpha_i=1+\left\lfloor\frac{i}{q}\right\rfloor$ is the smallest integer such that $q\cdot \alpha_i>i$, and $k=\left\lceil\frac{q(b-1)}{q+1}\right\rceil-1$ is the largest integer such that $\alpha_k\leq b-k-1$. So there are exactly
    $$2+\left\lceil\frac{q(b-1)}{q+1}\right\rceil-1=\left\lceil\frac{qb+1}{q+1}\right\rceil$$
    words covered by $0^a1^b$ in $\W_n^q$. Then we have
    $$    B_q(x)=\sum_{b=1}^{+\infty}\sum_{qa>b}^{+\infty}\left\lceil\frac{qb+1}{q+1}\right\rceil x^{a+b}
    =\frac{1}{1-x}\sum_{b=1}^{+\infty}\left\lceil\frac{qb+1}{q+1}\right\rceil x^{1+b+\left\lfloor\frac{b}{q}\right\rfloor}.
    $$
    \item We directly have
    $$D_q(x)=\sum_{b=1}^{+\infty}\sum_{qa>b}^{+\infty}x^{a+b}=\frac{1}{1-x}\sum_{b=1}^{+\infty}x^{1+b+\left\lfloor\frac{b}{q}\right\rfloor}.
    $$
\end{itemize}
    
\end{proof}

In the case where $q=c/d$ is rational, Proposition \ref{nb coverings} and Proposition \ref{closed form transfo} provide a closed form for $A_q(x),B_q(x)$ and $D_q(x)$. 

\begin{proposition}\label{closed form A B C}
    For $q=c/d$ a rational, we have
    $$A_q(x)=1+\frac{x}{1-x}+\frac{cx^{c+d}+(1-x)\sum_{r=0}^{c+d-1}\left\lfloor\frac{cr-1}{c+d}\right\rfloor x^r}{(1-x)(1-x^{c+d})}, \quad D_q(x)=\frac{\sum_{k=1}^c x^{1+\left\lfloor\frac{(c+d)k}{c}\right\rfloor}}{(1-x)(1-x^{c+d})},$$
    $$\mbox{and } B_q(x)=\frac{a_0x+a_1x^{1+\left\lfloor\frac{c+d}{c}\right\rfloor}+\ldots+a_{m}x^{1+\left\lfloor\frac{(c+d)m}{c}\right\rfloor}}{(1-x)(1-x^{c+d})(1-x^{(c+d)^2})},$$
    where $a_0,\ldots,a_m$ are the coefficients of the polynomial
    $$\left(x(1-x^{c+d})+cx^{c+d}+(1-x)\sum_{r=0}^{c+d-1}\left\lfloor\frac{cr+d-1}{c+d}\right\rfloor x^r\right)(1+x+\ldots+x^{c-1})(1+x^{c+d}+\ldots+x^{(c+d)(c-1)}).$$
\end{proposition}
\begin{proof}
It follows from Proposition \ref{nb coverings} and Fact \ref{fact partie entière} that $$A_q(x)=\sum_{n=1}^{+\infty}\left\lceil\frac{cn}{c+d}\right\rceil x^n=\sum_{n=1}^{+\infty}\left(1+\left\lfloor\frac{cn-1}{c+d}\right\rfloor\right) x^n=\frac{x}{1-x}+\sum_{n=0}^{+\infty}\left\lfloor\frac{cn-1}{c+d}\right\rfloor x^n-(-1).$$
The desired expression follows from Lemma \ref{gf (an/b)}. Using Proposition \ref{closed form transfo}, we directly deduce
$$D_q(x)=\frac{\sum_{k=1}^c x^{1+\left\lfloor\frac{(c+d)k}{c}\right\rfloor}}{(1-x)(1-x^{c+d})}.$$
Once again using Fact \ref{fact partie entière}, we have 
$$\left\lceil\frac{qb+1}{q+1}\right\rceil=\left\lceil\frac{cb+d}{c+d}\right\rceil=1+\left\lfloor\frac{cb+d-1}{c+d}\right\rfloor.$$
By Lemma \ref{gf (an/b)} we have, 
\begin{align*}
\sum_{b=1}^{+\infty}\left(1+\left\lfloor\frac{cb+d-1}{c+d}\right\rfloor\right)x^b&=\frac{x}{1-x}+\frac{cx^{c+d}+(1-x)\sum_{r=0}^{c+d-1}\left\lfloor\frac{cr+d-1}{c+d}\right\rfloor x^r}{(1-x)(1-x^{c+d})}\\
&=\frac{x(1-x^{c+d})+cx^{c+d}+(1-x)\sum_{r=0}^{c+d-1}\left\lfloor\frac{cr+d-1}{c+d}\right\rfloor x^r}{(1-x)(1-x^{c+d})}.
\end{align*}
Now it suffices to apply Proposition \ref{closed form transfo} to the above expression to obtain a closed form of $B_q(x)$. Observing that $\prod_{k=0}^{c-1}(1-e^{\frac{2i\pi k}{c}}x)(1-e^{\frac{2i\pi(c+d)k}{c}}x^{c+d})=(1-x^c)(1-x^{c(c+d)})$, and 
\begin{align*}
    \prod_{k=1}^{c-1}(1-e^{\frac{2i\pi k}{c}}x)(1-e^{\frac{2i\pi(c+d)k}{c}}x^{c+d})&=\frac{(1-x^c)(1-x^{c(c+d)})}{(1-x)(1-x^{c+d})}\\
    &=(1+x+\ldots+x^{c-1})(1+x^{c+d}+\ldots+x^{(c+d)(c-1)}),
\end{align*}
the expression of $B_q(x)$ then follows from Proposition \ref{closed form transfo}. 
\end{proof}

\begin{example} For instance, whenever $q\in\{1,2,1/2,1/3\}$, we obtain:
    \begin{align*}
        C_1(x)&=\frac{x(x^6 + x^3 - x^2 + x - 1)}{(x-1)(x^2 + 1)(x^2 + x - 1)^2}\\
        &=x+2x^2+5x^3+10x^4+20x^5+38x^6+70x^7+127x^8+228x^9+O(x^{10});\\[1em]
        C_2(x)&=\frac{x(x^{12} + x^{10} + 2x^8 + 2x^5 - x^4 - x^3 + 2x^2 - x - 1)}{(x-1)(x^6 + x^3 + 1)(x^3 + x^2 + x - 1)^2}\\
        &=x+4x^2+9x^3+22x^4+50x^5+108x^6+229x^7+476x^8+976x^9+O(x^{10});\\[1em]
        C_{1/2}(x)&=\frac{x(x^{12} + x^9 + x^7 - x^6 + x^4 - x^3 + x - 1)}{(x - 1)(x^6 + x^3 + 1)(x^3 + x - 1)^2}\\
        &=x+2x^2+3x^3+6x^4+11x^5+18x^6+30x^7+50x^8+81x^9+O(x^{10});\\[1em]
        C_{1/3}(x)&=\frac{x(x^{20} + x^{16} + x^{13} + x^9 - x^8 + x^5 - x^4 + x - 1)}{(x-1)(x^4 + 1)(x^8 + 1)(x^4 + x - 1)^2}\\
        &=x+2x^2+3x^3+4x^4+7x^5+12x^6+19x^7+28x^8+42x^9+O(x^{10}).
    \end{align*}
\end{example}

% Notice that for $q$ rational the generating  function $C_q(x)$ is also rational. 
%  Consider the following partial fraction decomposition for $q=1$, 
%  $$C_1(x)=1+\frac{-37 x +29}{25 \left(x^{2}+x -1\right)}+\frac{-2 x +4}{5 \left(x^{2}+x -1\right)^{2}}+\frac{1}{2 x -2}-\frac{x +7}{50 \left(x^{2}+1\right)}
% ,$$ 
% the coefficient $c_n$ of $x^n$, $n\geq 1$, in the  series expansion of $C_1(x)$ is given by 
% $$\frac{1}{25}\cdot\Big(\left(-4 n +33\right) F\! \left(n -1\right)+ \left(6 n -13\right)F \! \left(n \right)+4\left(n +1\right) F \! \left(n +1\right)\Big)-\frac{1}{2} - \Delta_n,$$
% where $F(n)$ is the $n$th term of the Fibonacci sequence, and $\Delta_n$ satisfies
% $$\Delta_n=\left\{\begin{array}{ll}
% -\frac{7}{50}(-1)^{n/2} & \mbox{ if } n \mbox{ is even,}\\[6pt]
% -\frac{1}{50}(-1)^{(n-1)/2}& \mbox{ otherwise.}
% \end{array}\right.$$

\begin{example}
   When $q$ is irrational, we do not have closed forms for the generating functions $A_q(x), B_q(x)$ and $D_q(x)$. However we can compute the first terms to deduce the first values in the Taylor expansion of $C_q(x)$. For instance, for $q\in\{\sqrt{2},\pi/4,e\}$, we obtain
    \begin{align*}
        C_{\sqrt{2}}(x)&=x+4x^2+9x^3+22x^4+46x^5+100x^6+207x^7+425x^8+856x^9+O(x^{10}),\\
        C_{\pi/4}(x)&=x+2x^2+5x^3+10x^4+20x^5+38x^6+70x^7+127x^8+224x^9+O(x^{10}),\\
        C_{e}(x)&=x+4x^2+12x^3+28x^4+67x^5+154x^6+343x^7+749x^8+1615x^9+O(x^{10}).
    \end{align*}
\end{example}
\begin{corollary}
    Let $q>0$ be a positive real number. Recall that $\Phi(q)^{-1}$ is the smallest root in modulus of the equation $D_q(x)=1$ (\cite{sergey2}). There exists a constant $c_q>0$ such that
    $$[x^n]C_q(x)\underset{n\to\infty}{\sim}c_q\cdot n\cdot \Phi(q)^n.$$
\end{corollary}
\begin{proof}
    It follows from routine asymptotic analysis. Indeed, from Proposition \ref{nb coverings}, $A_q(x)$ and $B_q(x)$ clearly have a radius of convergence of at least 1. So their smallest singularities are greater than $\Phi(q)^{-1}$, which is smaller than 1. So by Proposition \ref{nb coverings general}, the smallest singularity of $C_q(x)$ is then $\Phi(q)^{-1}$. Furthermore, this singularity has multiplicity 2, hence the factor $n$ in the asymptotic.
\end{proof}

\section{Intervals}\label{sec:intervals}

In this section, we focus on the enumeration of intervals in $\W_n^q$ for positive rational numbers $q$. An interval $I=\{u\in\W_n^q | \ v\leq u\leq w  \}$ is denoted $[v,w]$, it will be called \textit{prime} whenever $w$ is prime, i.e. $w=0^a1^b$ with $q\cdot a>b\geq1$. 

Let $P_q(x)$ be the generating function for the number of prime intervals in $\W_n^q$, and $I_q(x)$ the one for the number of all intervals. The following proposition establishes that the computation of $P_q(x)$ suffices to deduce the enumeration of all intervals.

\begin{prop}
 For any real $q>0$, the generating function  for the number of intervals is given by 
    $$I_q(x)=\frac{W_q(x)}{(1-x)(1-P_q(x))}.$$
\end{prop}
\begin{proof}
    Let $[v,w]$ be an interval in $\W_n^q$, and consider the decomposition of $w$:     $$w=1^m0^{a_1}1^{b_1}\ldots0^{a_k}1^{b_k}0^\ell,$$
    where $0^{a_i}1^{b_i}$ are prime factors. We take the decomposition of $v$ that is chosen to be compatible with that of $w$:
    $$v=u_1v_1\ldots v_ku_0,$$
    with $|u_1|=m$, $|u_0|=\ell$, and for all $i\in[1,k]$, $|v_i|=a_i+b_i$.
    
    By the definition of the relation $\leq$, we have $u_1\leq 1^m$, $u_0\leq 0^\ell$, and for all $i\in[1,k]$, $v_i\leq 0^{a_i}1^{b_i}$, so $u_1\in\mathcal{W}_m^q$, $u_0=0^\ell$, and $[v_i,0^{a_i}1^{b_i}]$ is a prime interval for all $i$. Conversely, given $u\in\mathcal{W}_m^q$  and some prime intervals $[v_i,0^{a_i}1^{b_i}]$, then $[u v_1\ldots v_k 0^\ell,1^m0^{a_1}1^{b_1}\ldots0^{a_k}1^{b_k}0^\ell]$ is an interval. The generating function for the number of elements $u\in \mathcal{W}_m^q$ is $W_q(x)$, the one for the number of sequences $0^\ell$ is $\frac{1}{1-x}$, and the one for the number of sequences of prime intervals is $\frac{1}{1-P_q(x)}$. The result then follows.
\end{proof}

Since $W_q(x)$ is known for any $q>0$, we now focus on the computation of $P_q(x)$.

% \begin{prop}\label{P_q general}
%     For $q\in\R_+$, we have
%     $$P_q(x)=\frac{1}{1-x}\sum_{b=1}^{+\infty}\left(1+\sum_{k=0}^{b-1}|W_k^q|\right)x^{b+\left\lceil\frac{b+1}{q}\right\rceil}.$$
% \end{prop}
% \begin{proof}
%     Let $a,b$ be such that $q\cdot a>b\geq1$, and $w=0^a1^b$. Then $v\in\mathcal{W}_{a+b}^q$ satisfies $v\leq w$ if and only if 
%     $v=0^{a+b}$ or $v=0^{a+k}1\nu$ for some $k\in [0,b-1]$, and $\nu\in\mathcal{W}_{b-k-1}^q$. Consequently,
%     \begin{align*}
%         P_q(x)&=\sum_{a=1}^{+\infty}\sum_{b=1}^{\lfloor qa-1\rfloor}\left(1+\sum_{k=0}^{b-1}|W_{b-k-1}^q|\right)x^{a+b}   =\sum_{a=1}^{+\infty}\sum_{b=1}^{\lfloor qa-1\rfloor}\left(1+\sum_{k=0}^{b-1}|W_k^q|\right)x^{a+b}\\
%         &=\sum_{b=1}^{+\infty}\left[\left(1+\sum_{k=0}^{b-1}|W_k^q|\right)x^b\sum_{a=\left\lceil\frac{b+1}{q}\right\rceil}^{+\infty}x^a\right]       =\sum_{b=1}^{+\infty}\left(1+\sum_{k=0}^{b-1}|W_k^q|\right)\frac{x^{b+\left\lceil\frac{b+1}{q}\right\rceil}}{1-x}.
%     \end{align*}
% \end{proof}

% We now give a formula to compute $P_q(x)$ when $q$ is rational. Recall that in this case we have an explicit expression for $W_q(x)$, see (\ref{eq:w_q rational}). With Lemma \ref{lemma periodicity} at hand, we can easily deduce the following formula for $P_q(x)$ when $q$ is rational.

\begin{thm}\label{p_q rational}
    Let $q=c/d$ be a positive rational number. Then we have
    % $$P_q(x)=\frac{1}{c(1-x)}\sum_{j=0}^{c-1}\sum_{i=0}^{c-1}\rho^{-ij}\Gamma_q(\rho^ix^{1+\frac{1}{q}})x^{1+\left\lfloor\frac{jd}{c}\right\rfloor-\frac{jd}{c}},$$
    % where $\rho=e^{2i\pi/c}$, and $\Gamma_q(x)=\frac{x(1+W_q(x))}{1-x}$.
    $$P_q(x)=\frac{\Gamma_{c/d}(x)}{1-x},$$
    where $\Gamma_{c/d}(x)=\sum_{n=0}^{+\infty} a_n x^{1+\left\lfloor\frac{(c+d)n}{c}\right\rfloor}$, with $\Gamma(x)=\frac{x(1+W_q(x))}{1-x}=\sum_{n=0}^{+\infty} a_n x^n$.
\end{thm}
\begin{proof}

Let $a,b$ be such that $q\cdot a>b\geq1$, and $w=0^a1^b$. Then $v\in\mathcal{W}_{a+b}^q$ satisfies $v\leq w$ if and only if $v=0^{a+b}$ or $v=0^{a+k}1\nu$ for some $k\in [0,b-1]$, and $\nu\in\mathcal{W}_{b-k-1}^q$. Consequently,
    \begin{align*}
        P_q(x)&=\sum_{b=1}^{+\infty}\sum_{a\geq 1+\left\lfloor\frac{bd}{c}\right\rfloor}\left(1+\sum_{k=0}^{b-1}|\mathcal{W}_{b-k-1}^q|\right)x^{a+b}  
        =\sum_{b=1}^{+\infty}\sum_{a\geq  1+\left\lfloor\frac{bd}{c}\right\rfloor}\left(1+\sum_{k=0}^{b-1}|\mathcal{W}_k^q|\right)x^{a+b}\\
        &=\sum_{b=1}^{+\infty}\left[\left(1+\sum_{k=0}^{b-1}|\mathcal{W}_k^q|\right)x^b\sum_{a= 1+\left\lfloor\frac{bd}{c}\right\rfloor}^{+\infty}x^a\right]       =\sum_{b=1}^{+\infty}\left(1+\sum_{k=0}^{b-1}|\mathcal{W}_k^q|\right)\frac{x^{1+b+\left\lfloor\frac{bd}{c}\right\rfloor}}{1-x}.
    \end{align*}
Since $\sum_{b=1}^{+\infty}\left(1+\sum_{k=0}^{b-1}|\mathcal{W}_k^q|\right)x^b=\Gamma(x)$, the result directly follows from the definition of $\Gamma_{c/d}(x)$ (see Lemma \ref{lemma periodicity}).    
\end{proof}

\begin{corollary}
    Let $q=c/d$ be a positive rational number, and define the polynomial $\Pi_q~=~1-x^{c+d}-\sum_{i=0}^{c-1}x^{1+i+\left\lfloor\frac{i}{q}\right\rfloor}$.
    \begin{enumerate}
        \item Let $\rho=e^{\frac{2i\pi}{c}}$ and $m=c(c+d+2)$. Then $$P_q(x)=\frac{\sum_{k=0}^m a_kx^{1+k+\left\lfloor\frac{k}{q}\right\rfloor}}{(1-x)(1-x^{c+d})^2\prod_{k=0}^{c-1}\Pi_q(\rho^kx^{1+\frac{d}{c}})},$$
        where $\sum_{k=0}^m a_kx^k=x\left(1-x^{c+d}+(1-x)\Pi_q(x)\right)(1+x+\ldots+x^{c-1})^2\prod_{k=1}^{c-1}\Pi_q(\rho^kx)$.
        \item Let $r=e^{\frac{2i\pi}{c+d}}$. There exists a nonzero polynomial $T\in\mathbb{C}[X]$ with degree at most $c+d-1$ such that
    $$[x^n]P_q(x)\underset{n\to\infty}{\sim}T(r^n)\cdot\Phi(q)^{\frac{cn}{c+d}}.$$
    \end{enumerate}
\end{corollary}
\begin{proof}
    It directly follows from Theorem \ref{p_q rational} and Proposition \ref{closed form transfo}, since $$\Gamma(x)=\frac{x(1-x^{c+d}+(1-x)\Pi_q(x))}{(1-x)^2\Pi_q(x)},$$
    and $[x^n]\Gamma(x)\underset{n\to\infty}{\sim}c_q\cdot \Phi(q)^n$ for some constant $c_q>0$.
\end{proof}

% \begin{corollary}
%     When $q$ is a positive integer, we have
%     $$P_q(x)=\frac{1}{q}\sum_{i=0}^{q-1}\frac{\Gamma_q(\rho^ix^{1+\frac{1}{q}})}{\rho^{-i}x^{-\frac{1}{q}}-1},$$
%     where $\rho=e^{2i\pi/q}$, and $\Gamma_q(x)=\frac{x(1+W_q(x))}{1-x}$.
% \end{corollary}

% \begin{proof}
%     Since $\left\lceil\frac{j}{q}\right\rceil=0$ for every $j\in [0,q-1]$, we have
%     \begin{align*}
%         P_q(x)&=\frac{1}{q(1-x)}\sum_{i=0}^{q-1}\Gamma_q(\rho^ix^{1+\frac{1}{q}})\sum_{j=0}^{q-1}\rho^{-ij}x^{1-\frac{j}{q}}=\frac{x}{q(1-x)}\sum_{i=0}^{q-1}\Gamma_q(\rho^ix^{1+\frac{1}{q}})\sum_{j=0}^{q-1}\left(\rho^{-i}x^{-\frac{1}{q}}\right)^j\\
%         &=\frac{x}{q(1-x)}\sum_{i=0}^{q-1}\Gamma_q(\rho^ix^{1+\frac{1}{q}})\frac{1-(\rho^{-i}x^{-\frac{1}{q}})^q}{1-\rho^{-i}x^{-\frac{1}{q}}}=\frac{1}{q}\sum_{i=0}^{q-1}\frac{\Gamma_q(\rho^ix^{1+\frac{1}{q}})}{\rho^{-i}x^{-\frac{1}{q}}-1}.
%     \end{align*}
% \end{proof}

% \begin{remark}
%      When $q$ is a positive integer, we have the following expression for $W_q(x)$:
%     $$ W_q(x)=\frac{1-x^{q+1}}{1-2x+x^{q+2}}.$$
% \end{remark}

\begin{example}
    When $q=1$, we have $P_1(x)=\frac{x^3(x^4-2)}{(x-1)^2(x+1)(x^4+x^2-1)}$, and
\begin{align*}
  I_1(x)&=\frac{(1-x)(1+x)^2(x^4+x^2-1)}{(x^6-2x^2-x+1)(x^2+x-1)}\\
  &=1+3x+6x^2+13x^3+27x^4+56x^5+115x^6+234x^7+474x^8+955x^9+O(x^{10}).  
\end{align*}

Furthermore, we have
\begin{align*}
    I_2(x)&=\frac{(1 - x)(x^9 + x^6 + 3x^3 - 1)(x^2 + x + 1)^2}{(x^{12} - x^7 + 2x^6 - 4x^3 - 2x^2 - x + 1)(x^3 + x^2 + x - 1)}\\
    &=1 + 3x + 9x^2 + 22x^3 + 57x^4 + 145x^5 + 363x^6 + 909x^7 + 2261x^8 + 5608x^9 + O(x^{10}),\\
    I_{2/3}(x)&=\scalemath{0.6}{\frac{(1-x)(x^4 + x^3 + x^2 + x + 1)^2(x^{25} + 2x^{20} + 3x^{15} + 2x^{10} + x^5 - 1)}{(x^{30} - x^{26} + x^{25} - x^{23} - 2x^{21} + x^{20} - 2x^{18} - 3x^{16} - x^{15} - 4x^{13} - 3x^{11} - x^{10} - 3x^8 - 2x^6 - 2x^5 - 2x^3 - x + 1)(x^5 + x^3 + x - 1)}}\\
        &=1 + 3x + 6x^2 + 13x^3 + 27x^4 + 51x^5 + 100x^6 + 196x^7 + 376x^8 + 723x^9 + O(x^{10}),\\
        I_{3/2}(x)&=\scalemath{0.5}{\frac{(1 - x)(x^{25} + 4x^{20} + 9x^{15} + 10x^{10} + 4x^5 - 1)(x^4 + x^3 + x^2 + x + 1)^2}{(x^5 + x^4 + x^2 + x - 1)(x^{30} - x^{26} + 3x^{25} - x^{22} - 3x^{21} + 5x^{20} - x^{19} - 4x^{17} - 5x^{16} + x^{15} - 5x^{14} - 6x^{12} - 5x^{11} - 6x^{10} - 8x^9 - 5x^7 - 3x^6 - 5x^5 - 4x^4 - 2x^2 - x + 1)}}\\
        &=1 + 3x + 9x^2 + 22x^3 + 57x^4 + 136x^5 + 336x^6 + 811x^7 + 1966x^8 + 4721x^9 + O(x^{10}).
    \end{align*}

% \begin{align*}
%     I_3(x)&=\frac{(1 - x)(x^{16} + x^{12} + x^8 + 7x^4 - 1)(x + 1)^2(x^2 + 1)^2}{(x^4 + x^3 + x^2 + x - 1)(x^{20} - x^{14} - 2x^9 + 6x^8 + x^7 - 8x^4 - 4x^3 - 2x^2 - x + 1)}\\
%     &=1 + 3x + 9x^2 + 27x^3 + 72x^4 + 199x^5 + 547x^6 + 1492x^7 + 4043x^8 + 10963x^9 + O(x^{10}).
% \end{align*}

% \begin{align*}
%     I_5(x)&=\frac{(1-x)(x + 1)^2(x^2 + x + 1)^2(x^2 - x + 1)^2(x^{36} + x^{30} + x^{24} + x^{18} - 9x^{12} + 31x^6 - 1)}{(x^{42} - x^{34} - 2x^{27} + x^{26} - 4x^{20} + 4x^{19} - 10x^{18} - x^{17} - 8x^{13} + 40x^{12} + 12x^{11} + 4x^{10} + x^9 - 32x^6 - 16x^5 - 8x^4 - 4x^3 - 2x^2 - x + 1)(x^6 + x^5 + x^4 + x^3 + x^2 + x - 1)}\\
%     &=1 + 3x + 9x^2 + 22x^3 + 57x^4 + 145x^5 + 363x^6 + 909x^7 + 2261x^8 + 5608x^9 + O(x^{10}).
% \end{align*}

\end{example}

\section{Meet-irreducible elements in $\W_n^q$ for rational numbers $q>0$}\label{sec:meet}

In this section we start by giving a closed form for the generating function of meet-irreducible elements in $\W_n^q$ for rational numbers $q$ in $]0,1]$, which is the simplest case among positive rational numbers. Then we give a general method to compute the generating function for rational numbers $q$ greater than 1.

\subsection{Meet-irreducible elements in $\W_n^q$ for a rational number $0<q\leq 1$}  The following easy fact enables us to consider only meet-irreducible elements starting with 0.

    \begin{fact}\label{adding 1s meet irreducible} Let $m,n\geq 1$ be two integers. If $w\in\W_n^q$, then $w$ is meet-irreducible if and only if $1^mw\in\W_{n+m}^q$ is meet-irreducible.
    \end{fact}
    
\begin{thm}\label{mi 0<q<=1}
    Let $q=c/d$ be a rational number with $0<q\leq 1$. The generating function $M(x)$ for the number of meet-irreducible elements in $\W_n^q$ is
    $$M(x)=\frac{x^{2 + \lfloor d/c\rfloor}(x-1)A - (x - 1)A^2  + x^2}{(x - 1)\left((1-Ax^{2 + \lfloor d/c\rfloor} + Ax^{3 + \lfloor d/c\rfloor} + (x - 1)(x+A-A^2 )\right)},$$
    with $A=\frac{\sum_{i=1}^c x^{1+i+\left\lfloor \frac{id}{c}\right\rfloor}}{1-x^{c+d}}$.
\end{thm}
\begin{proof} Taking into account Fact~\ref{adding 1s meet irreducible}, we focus on meet-irreducible elements starting with 0.

    Furthermore, let us consider the meet-irreducible elements ending with 1. These  elements can be written as concatenations of factors of the form \(0^a1^b\):
    $$w=0^{a_1}1^{b_1}\ldots 0^{a_k}1^{b_k},$$
    with $q\cdot a_i>b_i\geq1$ (which is equivalent to $a_i\geq 1+\left\lfloor\frac{db_i}{c}\right\rfloor$). 
    
    Note that each factor of such a meet-irreducible element satisfies $a_i\leq 1+\left\lfloor\frac{d(b_i+1)}{c}\right\rfloor$. Otherwise, it would mean that $0^{a-1}1^{b+1}\in\W_{a+b}^q$, and so $w$ would be smaller than both $1^{a_1+b_1}\ldots0^{a_i}1^{b_i}\ldots0^{a_k}1^{b_k}$, and $0^{a_1}1^{b_1}\ldots0^{a_i-1}1^{b_i+1}\ldots0^{a_k}1^{b_k}$, which are not comparable, contradicting the meet-irreducibility. 

    Thus, each factor in $w$ of the form $0^a1^b$ satisfies 
    \begin{equation}\label{eq:factor mi}
        1+\left\lfloor\frac{db}{c}\right\rfloor\leq a\leq 1+\left\lfloor\frac{d(b+1)}{c}\right\rfloor,
    \end{equation}
    and conversely, if $(a,b)$ satisfies (\ref{eq:factor mi}), then $0^a1^b$ is a prime factor of a meet-irreducible element. We now classify the factors satisfying (\ref{eq:factor mi}) into three sets $\mathtt{A},\mathtt{B},\mathtt{C}$ (not necessarily disjoint):
    \begin{itemize}
        \item $\mathtt{A}$ consists of the factors $0^a1^b$ such that        
        $$0^{a-1}1^b\not\in\W_{a+b-1}^q \mbox{ and } 0^a1^{b+1}\not\in\W_{a+b+1}^q,$$
        \item $\mathtt{B}$ consists of the factors $0^a1^b$ such that
        $$0^{a-1}1^b\in\W_{a+b-1}^q,$$
        \item $\mathtt{C}$ consists of the factors $0^a1^b$ such that
        $$0^a1^{b+1}\in\W_{a+b+1}^q.$$
    \end{itemize}

    For $0<q\leq1$, we have $\mathtt{C}\subseteq \mathtt{B}$ (see also (\ref{eq:caract B et C}) below), and $\mathtt{A}\cap \mathtt{B}=\emptyset$. The crucial condition on the meet-irreducible elements is that they cannot have a factor from $\mathtt{C}$ consecutively followed by a factor from $\mathtt{B}$ (what we call a \textit{pattern} $\mathtt{C}\mathtt{B}$). For instance, assume $w=0^{a_1}1^{b_1}0^{a_2}1^{b_2}$ with $0^{a_1}1^{b_1}\in \mathtt{C}$ and $0^{a_2}1^{b_2}\in \mathtt{B}$. Then $w$ is smaller than both $1^{a_1+b_1}0^{a_2}1^{b_2}$, and $0^{a_1}1^{b_1+1}0^{a_1-1}1^{b_1}$, which are incomparable, contradicting the meet-irreducibility. Similar examples can be built for words $w$ with any number of factors. Conversely, one can check that any word which is a product of factors $0^a1^b$ satisfying (\ref{eq:factor mi}), and avoiding the consecutive pattern $\mathtt{C}\mathtt{B}$, is meet-irreducible. Let $M_1(x)$ be the generating function of such words, i.e. meet-irreducible words starting with 0 and ending with 1. Counting these words is equivalent to count non-empty sequences of factors from the disjoint sets $\mathtt{A},\mathtt{B}\backslash\mathtt{C},\mathtt{C}$ avoiding consecutive factors $uv$, where $u\in\mathtt{C}$, and $v\in(\mathtt{B}\backslash\mathtt{C})$ or $v\in\mathtt{C}$. It means that a factor from $\mathtt{C}$ is either at the end of the word or followed by a factor from $\mathtt{A}$. Then we have 
    $$M_1(x)=\frac{A+B+CA}{1-A-(B-C)-CA},$$
    where $A,B,C$ are the generating functions of the factors in $\mathtt{A},\mathtt{B},\mathtt{C}$, respectively.
    
    Now let us study meet-irreducible elements starting with 0 and ending with 0. They have the following form: 
    $$w=0^{a_1}1^{b_1}\ldots 0^{a_k}1^{b_k}0^\ell,$$
    with $1+\left\lfloor\frac{db_i}{c}\right\rfloor\leq a_i\leq 1+\left\lfloor\frac{d(b_i+1)}{c}\right\rfloor$, $w$ avoids the consecutive pattern $\mathtt{C}\mathtt{B}$, and $1\leq\ell\leq1+\left\lfloor\frac{d}{c}\right\rfloor$. Indeed, if $\ell>1+\left\lfloor\frac{d}{c}\right\rfloor$, we would have $w$ smaller than both $1^{a_1+b_1}\ldots 0^{a_k}1^{b_k}0^\ell$ and $0^{a_1}1^{b_1}\ldots 0^{a_k}1^{b_k}0^{\ell-1}1$, which are not comparable. Moreover, $0^{a_k}1^{b_k}\not\in \mathtt{C}$, otherwise $w$ would be smaller than both $1^{a_1+b_1}\ldots 0^{a_k}1^{b_k}0^\ell$, and $0^{a_1}1^{b_1}\ldots 0^{a_k}1^{b_k+1}0^{\ell-1}$, which are not comparable. Conversely, such words $w$ are indeed meet-irreducible. Since the generating function for words counted by $M_1(x)$ and not ending with $\mathtt{C}$ is $\frac{A+(B-C)+CA}{1-A-(B-C)-CA}$, we deduce that the generating function $M_0(x)$ for meet-irreducible elements starting with 0 and ending with 0 is 
    $$M_0(x)=x+\ldots+x^{1+\left\lfloor\frac{d}{c}\right\rfloor}+\frac{A+(B-C)+CA}{1-A-(B-C)-CA}(x+\ldots+x^{1+\left\lfloor\frac{d}{c}\right\rfloor}).$$
    Using Fact \ref{adding 1s meet irreducible}, we deduce that the generating function $M(x)$ for all meet-irreducible elements is $$M(x)=\frac{M_0(x)+M_1(x)}{1-x}.$$
    To complete the proof it remains to compute the generating functions $A,B$ and $C$. Let us first characterize the factors belonging to $\mathtt{A},\mathtt{B}$ and $\mathtt{C}$, respectively. By checking the inequality defining $\mathcal{W}^q$, we have
    \begin{equation}\label{eq:caract B et C}
        0^a1^b\in \mathtt{B} \Leftrightarrow a\geq 2+\left\lfloor\frac{db}{c}\right\rfloor, \mbox{ and } 0^a1^b\in \mathtt{C} \Leftrightarrow a= 1+\left\lfloor\frac{d(b+1)}{c}\right\rfloor.
    \end{equation}
    Since for every $b\geq1$, $1+\left\lfloor\frac{db}{c}\right\rfloor<1+\left\lfloor\frac{d(b+1)}{c}\right\rfloor$, we have
    \begin{equation}\label{eq:caract A}
        0^a1^b\in \mathtt{A} \Leftrightarrow a= 1+\left\lfloor\frac{db}{c}\right\rfloor.
    \end{equation}
    By Lemma \ref{factor mi each length}, for each $n\geq2+\left\lfloor\frac{d}{c}\right\rfloor$, there exists a unique factor $0^a1^b$ of length $n$ such that $1+\left\lfloor\frac{db}{c}\right\rfloor\leq a\leq 1+\left\lfloor\frac{d(b+1)}{c}\right\rfloor$.
    % Furthermore, this factor is unique. Indeed, suppose that $w=1^m0^{a_1}1^{b_1}\ldots0^a1^b\ldots0^{a_k}1^{b_k}0^\ell$ and $w'=1^m0^{a_1}1^{b_1}\ldots0^{a'}1^{b'}\ldots0^{a_k}1^{b_k}0^\ell$ are both meet-irreducible, with $a+b=a'+b'=n$, and $a<a'$. Then we would have $w'\leq w$ and $w'\leq 1^{m+a_1+b_1}\ldots0^{a'}1^{b'}\ldots0^{a_k}1^{b_k}0^\ell$, contradicting the meet-irreducibility of $w'$, since $w$ and $1^{m+a_1+b_1}\ldots0^{a'}1^{b'}\ldots0^{a_k}1^{b_k}0^\ell$ are not comparable.
    So the generating function of all factors satisfying (\ref{eq:factor mi}) is $\frac{x^{2+\left\lfloor\frac{d}{c}\right\rfloor}}{1-x}$. Since $\mathtt{A}$ and $\mathtt{B}$ are disjoint, and $\mathtt{A}\cup \mathtt{B}$ contains all factors $0^a1^b$ satisfying (\ref{eq:factor mi}), $A+B=\frac{x^{2+\left\lfloor\frac{d}{c}\right\rfloor}}{1-x}$. By (\ref{eq:caract B et C}) and  (\ref{eq:caract A}), we have $$A=\sum_{b=1}^{+\infty}x^{1+b+\left\lfloor\frac{db}{c}\right\rfloor} \mbox{ and } C=\sum_{b=1}^{+\infty}x^{1+b+\left\lfloor\frac{d(b+1)}{c}\right\rfloor}=\frac{A}{x}-x^{1+\left\lfloor\frac{d}{c}\right\rfloor}.$$ Using Proposition \ref{closed form transfo}, we deduce $A=\frac{\sum_{i=1}^{c}x^{1+i+\left\lfloor\frac{di}{c}\right\rfloor}}{1-x^{c+d}}$. The desired expression of $M(x)$ follows after plugging $B=\frac{x^{2+\left\lfloor\frac{d}{c}\right\rfloor}}{1-x}-A$ and $C=\frac{A}{x}-x^{1+\left\lfloor\frac{d}{c}\right\rfloor}$ in $M_0(x)$ and $M_1(x)$, and simplifying.
\end{proof}

\begin{remark}
    If $w$ is meet-irreducible and $w\lessdot v$ with $v\ne1^n$, then $v$ is meet-irreducible, which means that the set of meet-irreducible elements, with $1^n$, forms an upper ideal.
\end{remark}

\subsection{Meet-irreducible elements in $\W_n^q$ for a rational number $q>1$} The enumeration of meet-irreducible elements for $q>1$ proceeds similarly to the case $0<q\leq1$ discussed in Theorem \ref{mi 0<q<=1}, although it is slightly more involved. As seen in Theorem \ref{mi 0<q<=1}, the structure of meet-irreducible elements for $0<q\leq1$ can be described in terms of words over a three-letter alphabet avoiding a consecutive pattern of length 2. As we will see in this section, for $q>1$, the structure of meet-irreducible elements corresponds to that of words over an alphabet  of $2\lceil q\rceil+1$ letters avoiding $\lceil q\rceil^2+2\lceil q\rceil-1$ consecutive patterns of length 2.

\textbf{The general pattern avoided.}
Fact \ref{adding 1s meet irreducible} still holds, and we start by focusing on words starting with 0 and ending with 1. Still as in the proof of Theorem \ref{mi 0<q<=1}, such meet-irreducible words are products of factors $0^a1^b$, with $1+\left\lfloor\frac{b}{q}\right\rfloor\leq a\leq 1+\left\lfloor\frac{b+1}{q}\right\rfloor$. Among these factors, we distinguish the same 3 sets:
\begin{itemize}
        \item $\mathtt{A}$ consists of the factors $0^a1^b\not\in\mathtt{C}$, 
        \item $\mathtt{B}$ consists of the factors $0^a1^b$ such that
        $$0^{a-1}1^b\in\W_{a+b-1}^q, \mbox{ and } a\geq2,$$
        \item $\mathtt{C}$ consists of the factors $0^a1^b$ such that
        $$0^a1^{b+1}\in\W_{a+b+1}^q.$$
    \end{itemize}
Note that we added the condition $a\geq2$ in $\mathtt{B}$, which was always guaranteed for $q\leq1$, but no longer for $q>1$. Furthermore, since $q>1$, we have $\mathtt{B}\subseteq \mathtt{C}$, in contrast to  the case $0<q\leq 1$ (see (\ref{eq:caract B et C q>1}) below). As in the proof of Theorem \ref{mi 0<q<=1}, the crucial condition on meet-irreducible elements is that they cannot have a factor from $\mathtt{C}$ consecutively followed by a factor from $\mathtt{B}$.

\textbf{The additional patterns avoided.}
In the case $q>1$, the avoidance of the pattern $\mathtt{C}\mathtt{B}$ is no longer sufficient to characterize meet-irreducible elements starting with 0 and ending with 1. Actually, the factors having only one 0, namely $01,011,\ldots,01^{\lceil q\rceil-1}$, play a special role (there are no such factors for $q\leq1$).

% Let us change slightly our definition of $\mathtt{A},\mathtt{B}$ and $\mathtt{C}$:
% \begin{itemize}
%         \item $\mathtt{A}$ consists of the factors $0^a1^b\not\in\mathtt{C}$, 
%         \item $\mathtt{B}$ consists of the factors $0^a1^b$ such that
%         $$0^{a-1}1^b\in\W_{a+b-1}^q, \mbox{ and } a\geq2,$$
%         \item $\mathtt{C}$ consists of the factors $0^a1^b$ such that
%         $$0^a1^{b+1}\in\W_{a+b+1}^q.$$
%     \end{itemize}
From the definitions of $\mathtt{A},\mathtt{B}$ and $\mathtt{C}$, we see that $01,\ldots,01^{\lceil q\rceil-2}\in \mathtt{C}$, and $01^{\lceil q\rceil-1}\in \mathtt{A}$. These factors play a special role since $01^i$ cannot follow a factor $0^a1^b$ such that $0^a1^{b+i+1}$ is $q$-decreasing, which would create an extra upper cover. So let us define new sets $\mathtt{D}_i$ for $1\leq i\leq \lceil q\rceil-1$:
\begin{itemize}
    \item $\mathtt{D}_i$ consists of the factors $0^a1^b$ such that $$0^a1^{b+i+1}\in\W_{a+b+i+1}^q.$$  
\end{itemize}
Then the patterns $\mathtt{D}_i01^i$ (i.e. factors in $\mathtt{D}_i$ followed by $01^i$)  are also forbidden in a meet-irreducible element.

\textbf{Forming the letters.}
Now we have described all the patterns avoided by a meet-irreducible element. However, there is one more step before fully converting our problem into a pattern avoiding problem on words. Indeed, the sets we consider, $\mathtt{A}$, $\mathtt{B}$, $\mathtt{C}$, the $\mathtt{D}_i$'s and the $01^i$'s are not disjoint. So we want to split them into disjoint sets in order to express properly the patterns avoided. One can check the following inclusions:
$$\mathtt{D}_{\lceil q\rceil-1}\subseteq \mathtt{B}\subseteq \mathtt{D}_{\lceil q\rceil-2}\subseteq \mathtt{D}_{\lceil q\rceil-3}\subseteq\ldots\subseteq \mathtt{D}_2\subseteq \mathtt{D}_1\subseteq \mathtt{C}.$$
Moreover, $01^{\lceil q\rceil-1}\in \mathtt{A}$, $01^{\lceil q\rceil-2}\in \mathtt{C}\backslash \mathtt{D}_1$, and $01^{\lceil q\rceil-2-i}\in \mathtt{D}_i\backslash \mathtt{D}_{i+1}$ for $1\leq i\leq \lceil q\rceil-3$. We then define the following $2\lceil q\rceil+1$ disjoint sets that will correspond to our letters:
$$\begin{array}{ccc}
    \mathtt{a}=\mathtt{A}\backslash\{01^{\lceil q\rceil-1}\}, & \mathtt{b}=\{01^{\lceil q\rceil-1}\},  \\
    \mathtt{c}=\mathtt{C}\backslash (\mathtt{D}_1\cup\{01^{\lceil q\rceil-2}\}), & \mathtt{d}=\mathtt{B}\backslash \mathtt{D}_{\lceil q\rceil-1}, & \mathtt{e}=\mathtt{D}_{\lceil q\rceil-1}, \\
    \mathtt{f}_i=\{01^i\} & \mbox{ for } 1\leq i\leq \lceil q\rceil-2,\\
    \mathtt{g}_i=\mathtt{D}_i\backslash (\mathtt{D}_{i+1}\cup\{01^{\lceil q\rceil-2-i}\}) & \mbox{ for } 1\leq i\leq\lceil q\rceil-3, & \mbox{ and } \mathtt{g}_{\lceil q\rceil-2}=\mathtt{D}_{\lceil q\rceil-2}\backslash \mathtt{B}.
\end{array}$$
Let also $\aaa,\bbb, \ccc,\ddd,\eee,\fff_1,\ldots,\fff_{\lceil q\rceil-2},\gggg_1,\ldots,\gggg_{\lceil q\rceil-2}$ denote, respectively, the generating functions of the corresponding sets.

\textbf{The matrix of forbidden patterns.}
We now have $2\lceil q\rceil+1$ letters, $\mathtt{a}$, $\mathtt{b}$, $\mathtt{c}$, $\mathtt{d}$, $\mathtt{e}$, $\mathtt{f}_1,\ldots,\mathtt{f}_{\lceil q\rceil-2}$, $\mathtt{g}_1,\ldots,\mathtt{g}_{\lceil q\rceil-2}$. The avoidance of $\mathtt{C}\mathtt{B}$, and $\mathtt{D}_i01^i$ for $1\leq i\leq \lceil q\rceil-1$ then translates into the following $\lceil q\rceil^2+2\lceil q\rceil-1$ forbidden patterns in terms of letters: $$\mathtt{c}\mathtt{d}, \mathtt{c}\mathtt{e}, \mathtt{d}\mathtt{d}, \mathtt{d}\mathtt{e}, \mathtt{e}\mathtt{b}, \mathtt{e}\mathtt{d}, \mathtt{e}\mathtt{e},$$ 
$$\mathtt{d}\mathtt{f}_i,\, \mathtt{e}\mathtt{f}_i,\,  \mathtt{f}_i \mathtt{d},\,  \mathtt{f}_i \mathtt{e},\,  \mathtt{g}_i \mathtt{d},\,  \mathtt{g}_i \mathtt{e},\, \,  \mathtt{g}_i\mathtt{f}_i,\, \mathtt{g}_i\mathtt{f}_{i-1},\ldots,\mathtt{g}_i\mathtt{f}_1 \mbox{ for } 1\leq i\leq\lceil q\rceil-2,$$ 
$$ \mbox{ and } \ \mathtt{f}_i\mathtt{f}_{\lceil q\rceil-2-i},\mathtt{f}_i\mathtt{f}_{\lceil q\rceil-3-i},\ldots,\mathtt{f}_i\mathtt{f}_1 \ \mbox{ for } 1\leq i\leq\lceil q\rceil-3.$$

For $\ell\in\{\mathtt{a},\mathtt{b},\mathtt{c},\mathtt{d},\mathtt{e},\mathtt{f}_1,\ldots,\mathtt{f}_{\lceil q\rceil-2},\mathtt{g}_1,\ldots,\mathtt{g}_{\lceil q\rceil-2}\}$, let $M_\ell$ denote the generating function of words on these $2\lceil q\rceil+1$ letters, avoiding these $\lceil q\rceil^2+2\lceil q\rceil-1$ patterns, and ending with the letter $\ell$. The patterns then induce structure on the $M_\ell$'s, for instance the letter $\mathtt{a}$ can be preceded by any letter, so 
$$M_\mathtt{a}=\aaa+(M_\mathtt{a}+M_\mathtt{b}+M_\mathtt{c}+M_\mathtt{d}+M_\mathtt{e}+M_{\mathtt{f}_1}+\ldots+M_{\mathtt{f}_{\lceil q\rceil-2}}+M_{\mathtt{g}_1}+\ldots+M_{\mathtt{g}_{\lceil q\rceil-2}})\cdot \aaa.$$
The letter $\mathtt{d}$ can be preceded only by $\mathtt{a}$ or $\mathtt{b}$, so
$$M_\mathtt{d}=\ddd+(M_\mathtt{a}+M_\mathtt{b})\cdot\ddd.$$
By doing this for each letter, we obtain a system on the $M_\ell$'s. By solving it we obtain an expression of each $M_\ell$ in terms of $\aaa,\bbb,\ldots,\gggg_{\lceil q\rceil-2}$, that we compute in the next paragraph. The generating function of meet-irreducible elements starting with 0 and ending with 1 is then
$$M_1=M_\mathtt{a}+M_\mathtt{b}+M_\mathtt{c}+M_\mathtt{d}+M_\mathtt{e}+M_{\mathtt{f}_1}+\ldots+M_{\mathtt{f}_{\lceil q\rceil-2}}+M_{\mathtt{g}_1}+\ldots+M_{\mathtt{g}_{\lceil q\rceil-2}}.$$
A meet-irreducible elements starts with 0 and ends with 0 if and only if it is 0 or it has the form $\ell_1\ldots\ell_k0$ with $k\geq1$, $\ell_1\ldots\ell_k$ avoids the above patterns, and $\ell_k\subseteq \mathtt{A}$, i.e. $\ell_k\in\{\mathtt{a},\mathtt{b}\}$. The generating function of meet-irreducible elements starting with 0 and ending with 0 is then
$$M_0=(1+M_\mathtt{a}+M_\mathtt{b})x.$$
By Fact \ref{adding 1s meet irreducible}, the generating function for all meet-irreducible elements is
$$M=\frac{M_0+M_1}{1-x}.$$

\textbf{Computing the generating functions}. To conclude, we compute the generating functions $A,B,C,D_i$ of the sets $\mathtt{A},\mathtt{B},\mathtt{C},\mathtt{D}_i$ in the next theorem.
\begin{thm}
    Let $q=c/d$ be a positive rational number, with $c$ and $d$ relatively prime. For $1\leq i\leq d$, let $a_i\in\{1,\ldots,c\}$ be such that $a_i=d^{-1}(c-d-1+i)\bmod{c}$, and for $0\leq m\leq c-d-1$, let $b_m\in\{1,\ldots,c\}$ be such that $b_m=d^{-1}m-1\bmod{c}$. Then
    $$A(x)=\frac{\sum_{i=1}^d x^{1+a_i+\left\lfloor\frac{a_i}{q}\right\rfloor}}{1-x^{c+d}}, \quad B(x)=\frac{\sum_{i=1}^d x^{2+a_i+\left\lfloor\frac{a_i}{q}\right\rfloor}}{1-x^{c+d}},\quad C(x)=\frac{\sum_{i=2}^{c+1}x^{i+\left\lfloor\frac{i}{q}\right\rfloor}}{1-x^{c+d}},$$
    $$\mbox{ and } \ D_i(x)=\frac{\sum_{m=0}^{c-1-di}x^{1+b_m+\left\lfloor\frac{b_m+1}{q}\right\rfloor}}{1-x^{c+d}} \  \mbox{ for } 1\leq i\leq\lceil q\rceil -1.$$
\end{thm}
\begin{proof}
By checking the inequality defining $\mathcal{W}^q$, we have
    \begin{equation}\label{eq:caract B et C q>1}
        0^a1^b\in \mathtt{B} \Leftrightarrow a\geq 2+\left\lfloor\frac{b}{q}\right\rfloor, \mbox{ and } 0^a1^b\in \mathtt{C} \Leftrightarrow a= 1+\left\lfloor\frac{b+1}{q}\right\rfloor.
    \end{equation}
    Furthermore, when $q>1$, we have  either $\left\lfloor\frac{b+1}{q}\right\rfloor=\left\lfloor\frac{b}{q}\right\rfloor$ or $\left\lfloor\frac{b+1}{q}\right\rfloor=1+\left\lfloor\frac{b}{q}\right\rfloor$, and $\mathtt{B}\subseteq \mathtt{C}$. To characterize the factors belonging to $\mathtt{A}$, we need to be careful when $\left\lfloor\frac{b}{q}\right\rfloor=\left\lfloor\frac{b+1}{q}\right\rfloor$:
    \begin{equation}\label{eq:caract A q>1}
        0^a1^b\in \mathtt{A} \Leftrightarrow a= 1+\left\lfloor\frac{b}{q}\right\rfloor\mbox{ and } \left\lfloor\frac{b+1}{q}\right\rfloor=1+\left\lfloor\frac{b}{q}\right\rfloor.
    \end{equation}
Indeed, if $a= 1+\left\lfloor\frac{b}{q}\right\rfloor=1+\left\lfloor\frac{b+1}{q}\right\rfloor$, then $0^a1^b$ belongs to $\mathtt{C}$, and not to $\mathtt{A}$. Let us first compute $C$: by (\ref{eq:caract B et C q>1}) and Proposition \ref{closed form transfo}
$$C(x)=\sum_{b=1}^{+\infty}x^{1+b+\left\lfloor\frac{b+1}{q}\right\rfloor}=\sum_{b=2}^{+\infty}x^{b+\left\lfloor\frac{b}{q}\right\rfloor}=\frac{\sum_{i=2}^{c+1}x^{i+\left\lfloor\frac{i}{q}\right\rfloor}}{1-x^{c+d}}.$$
By Lemma \ref{factor mi each length}, the generating function for all factors $0^a1^b$ such that $1+\left\lfloor\frac{b}{q}\right\rfloor\leq a\leq 1+\left\lfloor\frac{b+1}{q}\right\rfloor$ is $\frac{x^2}{1-x}$. Since $\mathtt{B}\subseteq \mathtt{C}$ because $q>1$, we have by definition $A(x)=\frac{x^2}{1-x}-C(x)$ (below, as a bonus and a warm up for the $D_i$'s, we compute a closed form of $A$ in another way, using a little bit of arithmetic). From (\ref{eq:caract A q>1}) and Fact \ref{carac b/q=(b+1)/q}, we deduce that 
$$A(x)=\sum_{\substack{b=1\\\left\lfloor\frac{b+1}{q}\right\rfloor=1+\left\lfloor\frac{b}{q}\right\rfloor}}^{+\infty}x^{1+b+\left\lfloor\frac{b}{q}\right\rfloor}=\sum_{\substack{b=1\\db\bmod{c}\geq c-d}}^{+\infty}x^{1+b+\left\lfloor\frac{b}{q}\right\rfloor}.$$
Since $c$ and $d$ are relatively prime, $d$ is invertible modulo $c$, so $A$ can be rewritten
$$A(x)=\sum_{b\in I}x^{1+b+\left\lfloor\frac{b}{q}\right\rfloor},$$
where $I=\{b\geq1 \ | \ b\bmod{c}\in\{d^{-1}k \ | \ k=c-d,c-d+1,\ldots,c-1\}\}$. 

For $1\leq i\leq d$, let $a_i\in\{1,\ldots,c\}$ be such that $a_i=d^{-1}(c-d-1+i)\bmod{c}$. Then $$\sum_{b\in I}x^b=\frac{x^{a_1}+\ldots+x^{a_d}}{1-x^c},$$ so by Proposition \ref{closed form transfo},
$$A(x)=\frac{\sum_{i=1}^d x^{1+a_i+\left\lfloor\frac{a_i}{q}\right\rfloor}}{1-x^{c+d}}.$$
Finally, using (\ref{eq:caract B et C q>1}) and $\left\lfloor\frac{b+1}{q}\right\rfloor\in\left\{\left\lfloor\frac{b}{q}\right\rfloor,1+\left\lfloor\frac{b}{q}\right\rfloor\right\}$, we actually have 
$$0^a1^b\in \mathtt{B} \Leftrightarrow \left\lfloor\frac{b+1}{q}\right\rfloor=1+\left\lfloor\frac{b}{q}\right\rfloor \mbox{ and } a= 2+\left\lfloor\frac{b}{q}\right\rfloor,$$
so that $$B(x)=\sum_{\substack{b=1\\\left\lfloor\frac{b+1}{q}\right\rfloor=1+\left\lfloor\frac{b}{q}\right\rfloor}}^{+\infty}x^{2+b+\left\lfloor\frac{b}{q}\right\rfloor}=x\cdot A(x)=\frac{\sum_{i=1}^d x^{2+a_i+\left\lfloor\frac{a_i}{q}\right\rfloor}}{1-x^{c+d}}.$$
Now it remains to compute the $D_i$'s. By checking the definition of $\mathcal{W}^q$, we have 
$$0^a1^b\in \mathtt{D}_i\Leftrightarrow a\geq1+\left\lfloor\frac{b+i+1}{q}\right\rfloor.$$
Since $1+\left\lfloor\frac{b}{q}\right\rfloor\leq a\leq 1+\left\lfloor\frac{b+1}{q}\right\rfloor$, we actually have 
$$0^a1^b\in \mathtt{D}_i\Leftrightarrow \left\lfloor\frac{b+i+1}{q}\right\rfloor=\left\lfloor\frac{b+1}{q}\right\rfloor \mbox{ and } a=1+\left\lfloor\frac{b+1}{q}\right\rfloor.$$
As a direct generalization of Fact \ref{carac b/q=(b+1)/q}, we have 
\begin{align*}
    \left\lfloor\frac{b+i+1}{q}\right\rfloor=\left\lfloor\frac{b+1}{q}\right\rfloor&\Leftrightarrow db=m-d\bmod{c}, \mbox{ for } m\in\{0,1,\ldots,c-1-di\}\\
    &\Leftrightarrow b=d^{-1}m-1\bmod{c}, \mbox{ for } m\in\{0,1,\ldots,c-1-di\}.
\end{align*}
We deduce that 
$$D_i(x)=\sum_{b\in J_i}x^{1+b+\left\lfloor\frac{b+1}{q}\right\rfloor},$$
where $J_i=\{b\geq1 \ | \ b\bmod{c}\in\{d^{-1}m-1 \ | \ m=0,1,\ldots,c-1-di\}\}$. 

For $0\leq m\leq c-d-1$, let $b_m\in\{1,\ldots,c\}$ be such that $b_m=d^{-1}m-1\bmod{c}$. Then for $1\leq i\leq \lceil q\rceil -1$,
$$\sum_{b\in J_i}x^b=\frac{x^{b_0}+\ldots+x^{b_{c-1-di}}}{1-x^c}.$$
Using Proposition \ref{closed form transfo}, we finally have 
$$D_i(x)=\frac{\sum_{m=0}^{c-1-di}x^{1+b_m+\left\lfloor\frac{b_m+1}{q}\right\rfloor}}{1-x^{c+d}}.$$
\end{proof}
\begin{example}
    Let us describe fully the case $q=5/2$. Here $\lceil q\rceil=3$, so there are 7 letters and 14 avoided patterns. From the previous computations, we have
    $$A(x)=\frac{x^3+x^6}{1-x^7},\quad B(x)=\frac{x^4+x^7}{1-x^7},\quad C(x)=\frac{x^2+x^4+x^5+x^7+x^8}{1-x^7},$$
    $$D_1(x)=\frac{x^4+x^7+x^8}{1-x^7},\quad D_2(x)=\frac{x^7}{1-x^7}.$$
    The corresponding generating functions for the letters are then
    $$\aaa=A-x^3,\quad \bbb=x^3,\quad \ccc=C-D_1-x^2,\quad \ddd=B-D_2,$$
    $$\eee=D_2,\quad \fff=x^2,\quad \gggg=D_1-B.$$
    The avoided patterns are $\mathtt{c}\mathtt{d}$, $\mathtt{c}\mathtt{e}$, $\mathtt{d}\mathtt{d}$, $\mathtt{d}\mathtt{e}$, $\mathtt{e}\mathtt{b}$, $\mathtt{e}\mathtt{d}$, $\mathtt{e}\mathtt{e}$, $\mathtt{d}\mathtt{f}$, $\mathtt{e}\mathtt{f}$, $\mathtt{f} \mathtt{d}$, $\mathtt{f} \mathtt{e}$, $\mathtt{g} \mathtt{d}$, $\mathtt{g} \mathtt{e}$, $\mathtt{g}\mathtt{f}$, which yield the following system:
    $$\left\{\begin{array}{cl}
       M_\mathtt{a} =& \aaa+(M_\mathtt{a}+M_\mathtt{b}+M_\mathtt{c}+M_\mathtt{d}+M_\mathtt{e}+M_\mathtt{f}+M_\mathtt{g})\cdot \aaa \\
       M_\mathtt{b} =& \bbb+(M_\mathtt{a}+M_\mathtt{b}+M_\mathtt{c}+M_\mathtt{d}+M_\mathtt{f}+M_\mathtt{g})\cdot \bbb \\
       M_\mathtt{c} =& \ccc+(M_\mathtt{a}+M_\mathtt{b}+M_\mathtt{c}+M_\mathtt{d}+M_\mathtt{e}+M_\mathtt{f}+M_\mathtt{g})\cdot \ccc \\
       M_\mathtt{d} =& \ddd+(M_\mathtt{a}+M_\mathtt{b})\cdot \ddd \\
       M_\mathtt{e} =& \eee+(M_\mathtt{a}+M_\mathtt{b})\cdot \eee \\
       M_\mathtt{f} =& \fff+(M_\mathtt{a}+M_\mathtt{b}+M_\mathtt{c}+M_\mathtt{f})\cdot \fff \\
       M_\mathtt{g} =& \gggg+(M_\mathtt{a}+M_\mathtt{b}+M_\mathtt{c}+M_\mathtt{d}+M_\mathtt{e}+M_\mathtt{f}+M_\mathtt{g})\cdot \gggg.
    \end{array}\right.$$
    Solving it, we obtain    
    \begin{align*}
    M_1&=M_\mathtt{a}+M_\mathtt{b}+M_\mathtt{c}+M_\mathtt{d}+M_\mathtt{e}+M_\mathtt{f}+M_\mathtt{g}\\
    &=\frac{-D_2 \left(B-D_1 \right) x^{5}+D_2 \left(B -C +1\right) x^{3}+\left(A B +D_1 \right) x^{2}-A B -A -C}{D_2 \left(B-D_1 \right) x^{5}-D_2 \left(B -C +1\right) x^{3}+\left(\left(1-A\right) B -D_1 \right) x^{2}+\left(A -1\right) B +A +C -1}\\
    &=-\frac{x^{2} \left(x^{2}+1\right) \left(x^{5}-x^{4}-x^{3}+x^{2}-1\right) \left(x^{11}+x^{10}-2 x^{8}-2 x^{7}+x^{5}-x^{4}+x^{2}+x +1\right)}{x^{20}-x^{18}-2 x^{17}+x^{15}+x^{14}+2 x^{12}+x^{10}+2 x^{9}-x^{8}-3 x^{7}-x^{6}-x^{5}-x^{3}-x^{2}+1},\\
    M_0&=(1+M_\mathtt{a}+M_\mathtt{b})x\\
    &=\frac{\left(\left(B-D_1 \right) x^{2}-B +C -1\right) x}{D_2 \left(B-D_1 \right) x^{5}-D_2 \left(B -C +1\right) x^{3}+\left(\left(1-A\right) B -D_1 \right) x^{2}+\left(A -1\right) B +A +C -1}\\
    &=-\frac{x \left(x^{7}-1\right) \left(x^{10}-x^{8}-x^{7}-x^{5}-x^{2}+1\right)}{x^{20}-x^{18}-2 x^{17}+x^{15}+x^{14}+2 x^{12}+x^{10}+2 x^{9}-x^{8}-3 x^{7}-x^{6}-x^{5}-x^{3}-x^{2}+1}.
    \end{align*}
    Finally, we deduce
    \begin{align*}
        M&=\frac{M_0+M_1}{1-x}\\
        &=\frac{x^{20}-2 x^{17}-2 x^{16}+x^{14}-2 x^{13}+2 x^{12}+4 x^{9}+x^{8}-2 x^{7}+x^{6}-x^{5}-x^{4}-x^{2}-x}{\left(x-1\right) \left(x^{20}-x^{18}-2 x^{17}+x^{15}+x^{14}+2 x^{12}+x^{10}+2 x^{9}-x^{8}-3 x^{7}-x^{6}-x^{5}-x^{3}-x^{2}+1\right)}\\
        &=x +2 x^{2}+3 x^{3}+6 x^{4}+9 x^{5}+13 x^{6}+23 x^{7}+34 x^{8}+52 x^{9}+O\! \left(x^{10}\right).
    \end{align*}
\end{example}

\section{Acknowledgement}
This research was funded, in part, by the Agence Nationale de la Recherche (ANR), grant ANR-22-CE48-0002 and by the Regional Council of Bourgogne-Franche-Comté.


\begin{thebibliography}{99}



\bibitem{Barc} E. Barcucci, A. Bernini, S. Bilotta
and R. Pinzani, Pattern avoiding and $q$-decreasing binary words,
\emph{RAIRO Theor. Inform. Appl.}, 59 (2025) 13.

\bibitem{arXivBKV} J.-L Baril, S. Kirgizov and V. Vajnovszki, Asymptotic bit frequency in {F}ibonacci words, {\em Pure Mathematics and Applications},  30(1) (2022), 23--30. 

\bibitem{bkv gray code} J.-L Baril, S. Kirgizov and V. Vajnovszki, Gray codes for Fibonacci $q$-decreasing words, \emph{Theoretical Computer Science}, 927 (2022), 120--132.




\bibitem{sergey2} S. Dovgal and S. Kirgizov, Structure and growth of $\R$-bonacci words,  {\em The Electronic Journal of Combinatorics}, 32(3) (2025), Article P3.32.

\bibitem{egecioglu irsic} \"{O}. E\u{g}ecio\u{g}lu and V. Ir\v{s}i\v{c}, Fibonacci-run graphs I: Basic properties, \emph{Discrete Applied Mathematics}, 295 (2021), 70--84.

\bibitem{Fein} M. Feinberg, Fibonacci–Tribonacci, \emph{Fibonacci Quart.}, 1 (1963) 71--74.

\bibitem{Grat} G.~Gr{\"a}tzer. {\em General Lattice Theory}, Second ed., Birkh{\"a}user, 1998.

% \bibitem{hsu} W.-J Hsu, Fibonacci cubes--a new interconnection Topology, \emph{IEEE Transactions on Parallel and Distributed Systems} 4 (1993), 3--12.

\bibitem{Serg} S. Kirgizov, $\Q$-bonacci words and numbers, \emph{Fibonacci Quarterly}, 60(5) (2022) 87–-195.

\bibitem{knuth} D.E. Knuth, The Art of Computer Programming: Sorting and Searching, 2nd edn., Vol. 3. \emph{Addison-Wesley Professional}, (1966).

%\bibitem{knuth42} D.E. Knuth, The Art of Computer Programming, Vol. 4, Fascicle 2, Generating All Tuples and Permutations, \emph{Addison-Wesley Professional}, (2005).



%\bibitem{l1} M. Lothaire, Combinatorics on words, \emph{Cambridge {U}niversity {P}ress}, 17 (1997).


\bibitem{Mile} E.P. Miles Jr, Generalized Fibonacci numbers and associated matrices, \emph{Amer. Math. Monthly}, 67 (1960) 745--752.


\bibitem{stanley} R.P. Stanley, Enumerative Combinatorics, Volume 1, \emph{Cambridge Studies in Advanced Mathematics}, (2011).

\bibitem{flaj}
P. Flajolet and R. Sedgewick, \emph{Analytic Combinatorics}, Cambridge University Press, 2005.


\bibitem{Wong} D. Wong, B. Liu, M. Im, Generating Cyclic 2-Gray Codes for Fibonacci $q$-Decreasing Words, \emph{WALCOM: Algorithms and Computation} (2024), 91--102.


\end{thebibliography}
\end{document}